\documentclass[12pt]{amsart}
\usepackage{amssymb}
\usepackage{amsthm}
\usepackage{amscd}
\usepackage{amsmath}
\usepackage{graphicx}

\theoremstyle{plain}
\newtheorem{lemma}{Lemma}[subsection]
\newtheorem{prop}[lemma]{Proposition}
\newtheorem{thm}[lemma]{Theorem}
\newtheorem{cor}[lemma]{Corollary}
\newtheorem{aplemma}{Lemma~A.\hspace{-1.5mm}}
\newtheorem{approp}{Proposition~A.\hspace{-1.5mm}}
\newtheorem{apthm}{Theorem~A.\hspace{-1.5mm}}
\newtheorem{apcor}{Corollary~A.\hspace{-1.5mm}}
\newtheorem{intthm}{Theorem}

\usepackage{textcomp}

\usepackage{pifont}

\usepackage{latexsym}

\usepackage[all]{xy}

\newtheorem{conj}[lemma]{Conjecture}

\theoremstyle{definition}

\newtheorem{rem}{Remark}
\newtheorem{rema}[lemma]{Remark}

\newtheorem{remb}{Remark}

\newtheorem{defi}[lemma]{Definition}
\newtheorem{exa}[lemma]{Example}
\newtheorem{aprem}{Remark~A.\hspace{-1.5mm}}
\newtheorem{apdefi}{Definition~A.\hspace{-1.5mm}}
\newcommand{\bde}{\begin{defi}}
\newcommand{\ede}{\end{defi}\vspace{1mm}}
\newcommand{\ble}{\begin{lemma}}
\newcommand{\ele}{\end{lemma}}
\newcommand{\bpr}{\begin{prop}}
\newcommand{\epr}{\end{prop}}
\newcommand{\bt}{\begin{thm}}
\newcommand{\et}{\end{thm}}
\newcommand{\bco}{\begin{cor}}
\newcommand{\eco}{\end{cor}}
\newcommand{\bre}{\begin{rem}}
\newcommand{\ere}{\end{rem}}
\newcommand{\brea}{\begin{rema}}
\newcommand{\erea}{\end{rema}\vspace{1mm}}
\newcommand{\breb}{\begin{remb}}
\newcommand{\ereb}{\end{remb}\vspace{1mm}}
\newcommand{\bex}{\begin{exa}}
\newcommand{\eex}{\end{exa}}
\newcommand{\bpf}{\begin{proof}}
\newcommand{\epf}{\end{proof}\vspace{1mm}}

\newcommand{\bade}{\begin{apdefi}}
\newcommand{\eade}{\end{apdefi}}
\newcommand{\bale}{\begin{aplemma}}
\newcommand{\eale}{\end{aplemma}}
\newcommand{\bapr}{\begin{approp}}
\newcommand{\eapr}{\end{approp}}
\newcommand{\bat}{\begin{apthm}}
\newcommand{\eat}{\end{apthm}}
\newcommand{\baco}{\begin{apcor}}
\newcommand{\eaco}{\end{apcor}}
\newcommand{\bare}{\begin{aprem}}
\newcommand{\eare}{\end{aprem}}

\newcommand{\be}{\begin{enumerate}}
\newcommand{\ee}{\end{enumerate}}
\newcommand{\bcd}{\[\begin{CD}}
\newcommand{\ecd}{\end{CD}\]}
\newcommand{\bit}{\begin{itemize}}
\newcommand{\eit}{\end{itemize}}
\newcommand{\bq}{\begin{quote}}
\newcommand{\eq}{\end{quote}}
\newcommand{\ba}{\begin{array}}
\newcommand{\ea}{\end{array}}
\newcommand{\mcA}{\mathcal{A}}
\newcommand{\mcB}{\mathcal{B}}
\newcommand{\mcC}{\mathcal{C}}
\newcommand{\mcD}{\mathcal{D}}
\newcommand{\mcE}{\mathcal{E}}
\newcommand{\mcF}{\mathcal{F}}
\newcommand{\mcG}{\mathcal{G}}
\newcommand{\mcH}{\mathcal{H}}

\newcommand{\mcL}{\mathcal{L}}

\newcommand{\mcO}{\mathcal{O}}

\newcommand{\mcT}{\mathcal{T}}

\newcommand{\mcV}{\mathcal{V}}

\newcommand{\mbC}{\mathbb{C}}

\newcommand{\mbF}{\mathbb{F}}

\newcommand{\mbN}{\mathbb{N}}

\newcommand{\mbS}{\mathbb{S}}

\newcommand{\mbU}{\mathbb{U}}

\newcommand{\mbZ}{\mathbb{Z}}
\newcommand{\mfA}{\mathfrak{A}}
\newcommand{\mfB}{\mathfrak{B}}
\newcommand{\mfC}{\mathfrak{C}}
\newcommand{\mfD}{\mathfrak{D}}
\newcommand{\mfE}{\mathfrak{E}}
\newcommand{\mfF}{\mathfrak{F}}
\newcommand{\mfG}{\mathfrak{G}}
\newcommand{\mfH}{\mathfrak{H}}

\newcommand{\mfL}{\mathfrak{L}}
\newcommand{\mfM}{\mathfrak{M}}

\newcommand{\mfO}{\mathfrak{O}}

\newcommand{\mfS}{\mathfrak{S}}

\newcommand{\mfX}{\mathfrak{X}}

\newcommand{\mfa}{\mathfrak{a}}

\newcommand{\mfc}{\mathfrak{c}}
\newcommand{\mfd}{\mathfrak{d}}
\newcommand{\mfe}{\mathfrak{e}}

\newcommand{\mfg}{\mathfrak{g}}
\newcommand{\mfh}{\mathfrak{h}}
\newcommand{\mfi}{\mathfrak{i}}

\newcommand{\mfk}{\mathfrak{k}}
\newcommand{\mfl}{\mathfrak{l}}

\newcommand{\mfn}{\mathfrak{n}}

\newcommand{\mfp}{\mathfrak{p}}

\newcommand{\mfs}{\mathfrak{s}}
\newcommand{\mft}{\mathfrak{t}}
\newcommand{\mfu}{\mathfrak{u}}

\newcommand{\migi}{\rightarrow}
\newcommand{\longmigi}{\longrightarrow}

\newcommand{\isom}{\stackrel{\sim}{\migi}}

\newcommand{\migiincl}{\hookrightarrow}

\newcommand{\migisurj}{\twoheadrightarrow}

\newcommand{\mr}{\mathrm}
\newcommand{\hidden}[1]{\,}

\newcommand{\BB}{\mcB}
\newcommand{\ZZZ}{{^{\mr{Zzz...}}}}

\begin{document}

\title{Duality for dormant opers}
\author{Yasuhiro Wakabayashi}
\date{}
\markboth{Yasuhiro Wakabayashi}{}
\maketitle
\footnotetext{Y. Wakabayashi: Graduate School of Mathematical Sciences, The University of Tokyo, 3-8-1 Komaba, Meguro, Tokyo,  153-8914, Japan;}
\footnotetext{e-mail: {\tt wkbysh@ms.u-tokyo.ac.jp};}
\begin{abstract}
In the present paper, we prove that on a fixed, pointed stable curve over a field   of characteristic $p> 0$,
there exists a {\it canonical} duality  between  dormant $\mathfrak{sl}_n$-opers ($1 < n <p-1$)  and   dormant $\mfs \mfl_{(p-n)}$-opers,
and  that there exists a unique (up to isomorphism) dormant $\mfs \mfl_{(p-1)}$-oper.

\end{abstract}
\tableofcontents 
\section*{Introduction}
The purpose of the present paper is to establish  a {\it canonical}  duality for dormant opers on a fixed algebraic curve  of characteristic $p>0$:
\vspace{2mm}
\begin{equation}
\text{\fbox{{\it dormant $\mfs \mfl_n$-opers}}} \ \   \rightleftharpoons  \ \ \text{\fbox{{\it dormant $\mfs \mfl_{(p-n)}$-opers}}}  \notag
\vspace{2mm}
\end{equation}
where $n$ is  an integer  with $1 < n < p-1$ and   $\mfs \mfl_n$ (resp., $\mfs \mfl_{(p-n)}$) denotes the special linear Lie algebra  of rank $n$ (resp., $p-n$).

\vspace{5mm}
\subsection*{0.1}\label{y01}

Recall that a  {\it dormant $\mfs \mfl_n$-oper}  is,  roughly speaking, a principal homogenous space over an algebraic curve  equipped with a connection satisfying certain conditions, including the condition that its  $p$-curvature vanishes identically.
Various properties of (dormant) $\mfs \mfl_n$-opers in characteristic $p>0$ and  $n =2$
  were first discussed
  by S. Mochizuki in ~\cite{Mzk1}.
If $n$ is general (but the underlying curve  is assumed to be unpointed and  smooth  over an algebraically closed field), then
the study of these objects
   has been carried out by K. Joshi, S. Ramanan, E. Z. Xia, J. K. Yu, C. Pauly, T. H. Chen, X. Zhu et al.
(cf. ~\cite{JRXY}, ~\cite{JP}, ~\cite{Jo14},  ~\cite{CZ}).
Also, formulations and background knowledge of dormant $\mfs \mfl_n$-opers (or, more generally, a dormant $(\mfg, \hslash)$-opers for a semisimple Lie algebra $\mfg$  and $\hslash \in k$) in the present  paper were discussed in the author's papers (cf. ~\cite{Wak5}, ~\cite{Wak6}).
As we explained briefly  in  ~\cite{Wak6}, \S\,0.2, 
dormant $\mfs \mfl_n$-opers and  their moduli, which are our principal objects,
contain diverse aspects and 
occur naturally in mathematics.
At any rate,  a detailed understanding  of them in   generalized setting will be of use  in  various  areas relevant to the theory of opers in positive characteristic.

\vspace{5mm}
\subsection*{0.2}\label{y02}

We shall describe  the main theorem of the present paper.
Let  $p$ a prime number, $n$ a positive integer with $n < p$,
 $(g,r)$ a pair of nonnegative integers satisfying the inequality $2g-2 +r >0$, $k$ a perfect field of characteristic $p$,  $S$ a $k$-scheme, and $\mfX_{/S} := (f : X \migi S, \{ \sigma_i : S \migi X \}_{i=1}^r)$  a pointed stable curve over $S$ of type $(g,r)$ (cf. \S\,\ref{z04}).
 Denote by $\mfc_n$ the GIT quotient of $\mfs \mfl_n$ by the adjoint action of $\mr{PGL}_n$ (= the projective linear group over $k$ of rank $n$).
 In \S\,\ref{z076},  we define  a certain subset $\mfc_n (\mbF_p)^\circledast$ (cf. (\ref{v0032})  for the precise definition of $\mfc_n (\mbF_p)^\circledast$)
  of 
 the set of $\mbF_p$-rational points $\mfc_n (\mbF_p)$ of $\mbF_p$. 
Let  $\vec{\rho}$ be an element of $(\mfc_n (\mbF_p)^\circledast)^{\times r}$ (i.e., the product of $r$ copies of $\mfc_n (\mbF_p)^\circledast$).
According to the discussion in \S\,\ref{z076},  
to each such  $\vec{\rho}$,  one may  associate  another element $\vec{\rho}^{\, \bigstar}$ of $(\mfc_n (\mbF_p)^\circledast)^{\times r}$.
(Here, if $r =0$, then we take $\vec{\rho} =\vec{\rho}^{\, \bigstar} = \emptyset$.)
We shall write 
\begin{equation}
 \mfO \mfp^\ZZZ_{\mfs \mfl_n, \mfX_{/S}} \ \ \left(\text{resp.,} \  \mfO \mfp^\ZZZ_{\mfs \mfl_n, \vec{\rho}, \mfX_{/S}}\right)
\end{equation}
for  the moduli stack classifying dormant $\mfs \mfl_n$-opers  (resp., dormant $\mfs \mfl_n$-opers  of radii $\vec{\rho}$) on $\mfX_{/S}$.
Then, the main theorem of the present paper is the  following assertion
(cf. Theorem \ref{z078}).

\vspace{3mm}
\begin{intthm} \label{ww033} \leavevmode\\
 \ \ \ Suppose that $n < p-1$.
 \begin{itemize}
\item[(i)]
 There exists a canonical isomorphism
 \begin{equation}
 \Theta_{\mfs \mfl_n, \mfX_{/S}}^\bigstar : \mfO \mfp^\ZZZ_{\mfs \mfl_n, \mfX_{/S}} \isom \mfO \mfp^\ZZZ_{\mfs \mfl_{(p-n)}, \mfX_{/S}} 
 \end{equation}
over $S$ satisfying that $\Theta_{\mfs \mfl_{(p-n)}, \mfX_{/S}}^\bigstar \circ   \Theta_{\mfs \mfl_n, \mfX_{/S}}^\bigstar = \mr{id}$.
\item[(ii)]
By restricting  $\Theta_{\mfs \mfl_n, \mfX_{/S}}^\bigstar$, we obtain  a canonical  isomorphism 
\begin{equation}
 \Theta_{\mfs \mfl_n, \vec{\rho}, \mfX_{/S}}^\bigstar : \mfO \mfp^\ZZZ_{\mfs \mfl_n, \vec{\rho}, \mfX_{/S}} \isom \mfO \mfp^\ZZZ_{\mfs \mfl_{(p-n)}, \vec{\rho}^{\, \bigstar}, \mfX_{/S}} 
\end{equation}
over $S$ satisfying that $\Theta_{\mfs \mfl_{(p-n)}, \vec{\rho},\mfX_{/S}}^\bigstar \circ   \Theta_{\mfs \mfl_n, \vec{\rho}^{\, \bigstar},\mfX_{/S}}^\bigstar = \mr{id}$.
\end{itemize}
\end{intthm}
\vspace{5mm}

Also, we obtain the following assertion (cf. Theorem \ref{z07jjj8}), which is a generalization (to the case of  pointed stable curves) of 
~\cite{Hoshi}, Theorem A (ii).
The unique dormant $\mfs \mfl_{(p-1)}$-oper asserted in {\it loc.\,cit.} was studied  explicitly by means of (the Frobenius pull-back of) the sheaf of locally exact $1$-forms. 

\vspace{3mm}
\begin{intthm} \label{ww03ggg3} \leavevmode\\
 \ \ \ The structure morphism $\mfO \mfp^\ZZZ_{\mfs \mfl_{(p-1)}, \mfX_{/S}} \migi S$ of $\mfO \mfp^\ZZZ_{\mfs \mfl_{(p-1)}, \mfX_{/S}}$ is an isomorphism.
That is to say, there exists a unique (up to isomorphism) dormant $\mfs \mfl_{(p-1)}$-oper on $\mfX_{/S}$.
\end{intthm}
\vspace{5mm}

One may apply Theorems A and B to the study toward explicit computations of 
the number of dormant $\mfs \mfl_n$-opers.
Indeed, let $\overline{\mfM}_{g,r}$ denote the moduli stack classifying pointed stable curves over $k$ of type $(g,r)$ and $\mfO \mfp^\ZZZ_{\mfs \mfl_n, \vec{\rho}, g,r}$ denote the moduli stack classifying pointed stable curves over $k$ of type $(g,r)$ equipped with a dormant $\mfs \mfl_n$-oper of radii $\vec{\rho}$ on it.
Then, Theorem A allows us to generalize the result of ~\cite{Wak5}, Theorem H, which give an explicit computation of the generic degree  of $\mfO \mfp^\ZZZ_{\mfs \mfl_n, \vec{\rho}, g,0}$ over $\overline{\mfM}_{g,0}$ (cf. Corollary \ref{z079}).
Here, recall that since $\mfO \mfp^\ZZZ_{\mfs \mfl_n, \vec{\rho}, g,0}$ is finite and generically \'{e}tale over $\overline{\mfM}_{g,r}$ (cf. ~\cite{Wak5}, Theorem G), this generic degree coincides with the number of dormant $\mfs \mfl_n$-opers  on a sufficiently general curve.

Moreover, 
by combining results in the present paper with  results in  $p$-adic Teichm\"{u}ller theory, we will  have (in \S\,\ref{z0762}) a rather explicit understanding of the 
case where $n = p-2$.
In particular, we will discuss (cf. Corollary \ref{z08022}) the structure of the fusion ring $ \mfF^\ZZZ_{p, \mfs \mfl_{(p-2)}}$ (cf. (\ref{d0d0d}))
associated with 
the function $N^\ZZZ_{p, \mfs \mfl_{(p-2)}, 0}$ (cf. (\ref{e0e0e})) assigning, to each data of radii $\vec{\rho} \in (\mfc_n (\mbF_p)^\circledast)^{\times r}$,  the generic degree of $\mfO \mfp^\ZZZ_{\mfs \mfl_{(p-2)}, \vec{\rho}, 0,r}/\overline{\mfM}_{0,r}$.

\vspace{5mm}
\hspace{-4mm}{\bf Acknowledgement} \leavevmode\\
 \ \ \ The author would like to express his sincere gratitude to Professors  Kirti Joshi and Shinichi Mochizuki for their inspiring works concerning dormant opers.
 Also, the author  would like to thank the referee for  reading carefully his  manuscript and giving him  some helpful comments and suggestions. 




\vspace{6mm}
\section{Preliminaries} \label{z01}\vspace{3mm}

Throughout the present paper, let us fix a prime $p$,  a perfect field $k$ of characteristic $p$ (hence $\mbF_p := \mbZ/p \mbZ \subseteq k$), and
a pair  of nonnegative integers  $(g, r)$ satisfying that $2g-2+r >0$.


\vspace{5mm}
\subsection{} \label{z04}
For a field $k'$ over $\mbF_p$,
  we shall denote by
 \begin{equation}
 2^{k'} \  \left(\text{resp.},  \ 2^{k'}_{\sharp(-)=n}\right)
 \end{equation}
  the set of subsets of $k'$ (resp., the set of subsets of $k'$ with cardinality $n$).
 Also, denote by 
 \begin{equation}
 \mbN^{k'} \  \left(\text{resp.},  \ \mbN^{k'}_{\sharp (-)= n}\right)
 \end{equation}
  the set of multisets over $k'$ (resp., the set of multisets over $k'$ with cardinality $n$).
(For the definition and various notations  concerning a {\it multiset}, we refer to ~\cite{SIYS}.)
In particuler,  $2^{k'} \subseteq \mbN^{k'}$ and $2^{k'}_{\sharp (-)=n} \subseteq \mbN^{k'}_{ \sharp (-)=n}$.

The symmetric group $\mfS_n$ of $n$ letters acts,   by permutation,  on the product $k'^{\times n}$  of $n$ copies of $k'$.
The quotient set $\mfS_n \backslash k'^{\times n}$ may be identified with $\mbN^{k'}_{ \sharp (-)=n}$ and, in particular, we have a natural surjection  
\begin{equation} \label{w101}
k'^{\times n} \migisurj \left(\mfS_n \backslash k'^{\times n} \cong\right) \ \mbN^{k'}_{ \sharp (-)=n}.
\end{equation}

Let $\tau_0 : = [\tau_{0, 1}, \cdots, \tau_{0, n}]$ be an element of $ \mbN^{k'}_{ \sharp (-)=n}$  and $a \in k'$.
Then, we shall write
\begin{align} \label{h990}
& \hspace{14mm} \tau_0^{+a} :=  [\tau_{0, 1} +a, \cdots, \tau_{0, n}+a]  \\
& \left(\text{resp.,} \ \tau_0^{-a}  :=  [\tau_{0, 1} -a, \cdots, \tau_{0, n}-a] \ \right) \in \mbN^{k'}_{ \sharp (-)=n}.\notag
\end{align}
If, moreover,  $\tau_0$  is  a subset of $\mbF_p$ (i.e., an element of $2_{ \sharp (-)=n}^{\mbF_p} \subseteq  \mbN^{k'}_{ \sharp (-)=n}$), then
we shall write 
\begin{equation} \label{w03}
\tau_0^\triangleright : = \mbF_p \setminus \tau_0 \ \left(\subseteq \mbF_p\right), 
\end{equation}
and  
\begin{equation} \label{w04}
\tau^\triangledown_0  :=[-\tau_{0, 1}, -\tau_{0, 2}, \cdots, -\tau_{0, n}]  \ \left(\subseteq \mbF_p\right).
\end{equation}

Next, let  $r$ be a positive  integer,  and $\vec{\tau} := (\tau_{i})_{i=1}^r$ an $r$-tuple 
of multisets over $k'$, i.e., an element of the product $(\mbN^{k'})^{\times r}$ of $r$-copies of $\mbN^{k'}$.
Then, for each  $\vec{a} := (a_i)_{i=1}^r \in k'^{\times r}$, we shall write 
\begin{equation} \label{z003}
\vec{\tau}^{\,+\vec{a}} := (\tau_i^{+ a_i})_{i=1}^r \ \ \left(\text{resp.,} \ \vec{\tau}^{\,-\vec{a}} := (\tau_i^{- a_i})_{i=1}^r\right).
\end{equation}
If, moreover,  $\vec{\tau}$ lies in $(2_{ \sharp (-)=n}^{\mbF_p})^{\times r}$, then
 we shall write
\begin{equation} \label{w034}
\vec{\tau}^{ \, \triangleright} := (\tau_i^\triangleright)_{i=1}^r, \ \ \ \vec{\tau}^{\triangledown} := (\tau_i^\triangledown)_{i=1}^r, \ \ \text{and} \  \  \vec{\tau}^{\, \bigstar} := ((\tau_i^\triangleright)^\triangledown)_{i=1}^r,
\end{equation}
all of which are elements of $(2_{ \sharp (-)=n}^{\mbF_p})^{\times r}$.
One verifies immediately the equalities
\begin{equation} \label{z02}
(\vec{\tau}^{\, \bigstar})^\bigstar=\vec{\tau}
\end{equation}
and 
\begin{equation}  \label{z02ee}
(\vec{\tau}^{\, +\vec{a}})^{\bigstar}=(\vec{\tau}^{\, \bigstar})^{-\vec{a}}.
\end{equation}


\vspace{5mm}
\subsection{} \label{z03}

Let $T$ be a scheme  over $k$ ($\supseteq \mbF_p$) and $f : Y \migi T$ a scheme over $T$.
Denote by $F_T : T \migi T$ (resp., $F_Y : Y \migi Y$) the absolute Frobenius morphism of $T$ (resp., $Y$).
The  {\it Frobenius twist of $Y$ over $T$} is, by definition,  the base-change $Y^{(1)}_T$ ($:= Y \times_{T, F_T} T$) of $f : Y\migi T$ via
 $F_T : T \migi T$.
Denote by $f^{(1)} : Y^{(1)}_T \migi T$ the structure morphism of the Frobenius twist of $Y$ over $T$. 
The {\it relative Frobenius morphism of $Y$ over $T$} is  the unique morphism $F_{Y/T} : Y \migi Y^{(1)}_T$ over $T$ that fits into a commutative diagram of the form 
\begin{equation} \begin{CD}
Y @> F_{Y/T} >> Y^{(1)}_T @> \mr{id}_Y \times F_T >> Y
\\
@V f VV @V f^{(1)} VV @V f VV
\\
T @> \mr{id}_T >> T @> F_T >> T.
\end{CD} \end{equation}
Here,  the upper  composite in this diagram coincides with $F_Y$
 and the right-hand square is, by the definition of $Y^{(1)}_T$, cartesian.

\vspace{5mm}
\subsection{} \label{z04}

Denote by $\overline{\mfM}_{g,r}$ the moduli stack of $r$-pointed stable curves (cf. ~\cite{Kn2}, Definition 1.1) over $k$ of genus $g$ 
 (i.e., of type $(g,r)$), and by $f_\mr{\mft \mfa \mfu} : \mfC_{g,r} \migi \overline{\mfM}_{g,r}$ the tautological curve, with its $r$ marked points $\mfs_1 , \cdots, \mfs_r : \overline{\mfM}_{g,r} \migi \mfC_{g,r}$.
Recall (cf. ~\cite{Kn2}, Corollary 2.6 and Theorem 2.7; ~\cite{DM}, \S\,5) that $\overline{\mfM}_{g,r}$ may be represented by a geometrically connected, proper, and smooth Deligne-Mumford stack over $k$ of dimension $3g-3+r$.
Also, recall (cf. ~\cite{KaFu}, Theorem 4.5) that $\overline{\mfM}_{g,r}$ has a natural log structure given by the divisor at infinity,
where we shall denote the resulting log stack  by $\overline{\mfM}_{g,r}^{\mr{log}}$.
Also, by taking the divisor which is the union of the $\mfs_i$'s and the pull-back of the divisor at infinity of $\overline{\mfM}_{g,r}$, we obtain a log structure on $\mfC_{g,r}$; we denote  the resulting log stack by $\mfC^{\mr{log}}_{g,r}$.
$f_\mr{\mft \mfa \mfu} : \mfC_{g,r} \migi \overline{\mfM}_{g,r}$ extends naturally to a morphism 
$f_\mr{\mft \mfa \mfu}^\mr{log} : \mfC^\mr{log}_{g,r} \migi \overline{\mfM}^\mr{log}_{g,r}$ of log stacks.

Next,  let $S$ be a scheme, or more generally, a stack  over $k$ and 
\begin{equation}
\label{X}
 \mfX_{/S} : =(f :X \migi S, \{ \sigma_i : S \migi X\}_{i=1}^r)\end{equation}
 a pointed stable curve over $S$ of type $(g,r)$, consisting of a (proper) semi-stable curve $f : X \migi S$ over $S$ of genus $g$ and $r$ marked points $\sigma_i : S \migi X$ ($i = 1, \cdots, r$).
$\mfX_{/S}$ determines its classifying morphism $s : S \migi \overline{\mfM}_{g,r}$ and an isomorphism $X \isom S \times_{s, \overline{\mfM}_{g,r}} \mfC_{g,r}$ over $S$.
By pulling-back the  log structures of $\overline{\mfM}^{\mr{log}}_{g,r}$ and $\mfC^{\mr{log}}_{g,r}$,
we obtain log structures on $S$ and $X$ respectively; we denote  the resulting log stacks by $S^{\mr{log}}$ and $X^{\mr{log}}$. 
The structure morphism $f : X \migi S$ extends to a morphism $f^\mr{log} : X^{\mr{log}} \migi S^{\mr{log}}$ of log stacks, which is log smooth (cf. ~\cite{KaKa}, \S\,3;  ~\cite{KaFu}, Theorem 2.6). 
Write $\mcT_{X^\mr{log}/S^\mr{log}}$ for  the sheaf  of logarithmic derivations of $X^\mr{log}$ over $S^\mr{log}$  and $\Omega_{X^\mr{log}/S^\mr{log}} := \mcT_{X^\mr{log}/S^\mr{log}}^\vee$  for
 its dual, i.e., 
 the sheaf   of logarithmic differentials of $X^\mr{log}$ over $S^\mr{log}$.
 Both $\mcT_{X^\mr{log}/S^\mr{log}}$ and  $\Omega_{X^\mr{log}/S^\mr{log}}$  are  line bundles on $X$.
For each $i = 1, \cdots, r$,
one may obtain  the {\it residue isomorphism}, which is, by definition, an isomorphism
\begin{equation} \label{triv}
 \sigma_i^{*}(\Omega_{X^\mr{log}/S^{\mr{log}}}) \isom \mcO_{S}
 \end{equation}
given by assigning $1 \in \mcO_{S}$ to  any local section of the form $\sigma_i^{*}(d\mr{log}(x)) \in \sigma_i^{*}(\Omega_{X^\mr{log}/S^{\mr{log}}})$  (for   a local function  $x$ defining 
the closed subscheme $\sigma_i : S \migi X$ of $X$).

If $t : T \migi S$ is an $S$-scheme, then
we shall use the notation ``$\mfX_{/T}$" for  indicating 
the base-change of $\mfX_{/S}$ via $t$, i.e., the pointed stable curve 
\begin{equation} \label{BC44} \mfX_{/T} := (X \times_{S, t} T/T, \{ \sigma_i \times_S \mr{id}_T : (T =) \ S \times_{S,  t}T \migi X \times_{S, t} T\}_{i=1}^r) \end{equation}
over $T$.

\vspace{5mm}
\subsection{} \label{z06}

Let 
$S$, $\mfX_{/S}$ be as above and
 $\mcV$   a    vector bundle (i.e., a locally free coherent $\mcO_X$-module)  on $X$.
By an {\it $S$-connection} (resp., {\it $S$-log connection}) on  $\mcV$,  we mean (cf. ~\cite{Wak5}, \S\,4.1) an
$f^{-1}(\mcO_S)$-linear morphism
\begin{equation}
  \nabla : \mcV \migi \Omega_{X/S} \otimes \mcV \ \left(\text{resp.,} \  \nabla : \mcV \migi \Omega_{X^\mr{log}/S^\mr{log}} \otimes \mcV\right) 
  \end{equation}
 satisfying the condition that
 \begin{equation}  \nabla (a \cdot m) =  d (a) \otimes m + a \cdot \nabla (m)\end{equation}
 for local sections $a \in \mcO_X$ and $m \in \mcV$, where $d$ denotes the universal  derivation $\mcO_X \migi \Omega_{X/S}$ ($\subseteq \Omega_{X^\mr{log}/S^\mr{log}}$). 

If $\nabla$ is an $S$-log connection on $\mcV$, then
we shall write
 $\mr{det}(\nabla)$  (resp., $\nabla^\vee$)
 for the $S$-log connection on the determinant $\mr{det}(\mcV)$ (resp., the dual $\mcV^\vee$) of $\mcV$ induced naturally by $\nabla$.
Also, for $m \geq 1$, we shall write $\nabla^{\otimes m}$ for the  $S$-log connection on the $m$-fold tensor product $\mcV^{\otimes m}$ of $\mcV$ induced by $\nabla$.
If, moreover, we are given a vector bundle $\mcV'$ on $X$  and  an $S$-log connection $\nabla'$ on $\mcV'$, then
we shall write $\nabla \otimes \nabla'$ for the $S$-log connection on the tensor product $\mcV \otimes \mcV'$ induced by $\nabla$ and $\nabla'$.

A {\it log integrable vector bundle} on $\mfX_{/S}$ (of rank $m \geq 1$) is a pair $\mfF := (\mcF, \nabla_\mcF)$ consisting of a vector bundle $\mcF$ on $X$ (of rank $m$) and an $S$-log connection $\nabla_\mcF$ on $\mcF$.
 If, moreover,  $\mcF$ is  of  rank $1$, then we shall refer to such an  $\mfF$ as  a {\it log integrable line bundle} on $\mfX_{/S}$.  

Let $\mfF := (\mcF, \nabla_\mcF)$ and $\mfG := (\mcG, \nabla_\mcG)$ be log integrable vector bundles on $\mfX_{/S}$.
An {\it isomorphism of log integrable vector bundles} 
 from $\mfF$ to $\mfG$ is an isomorphism $\mcF \isom \mcG$ of $\mcO_X$-modules that is compatible with the respective $S$-log connections $\nabla_\mcF$ and $\nabla_\mcG$.

\vspace{5mm}

\subsection{} \label{z07}

We recall the definition of the $p$-curvature of a logarithmic connection (cf., e.g., ~\cite{Wak}, \S\,3).
Let 
$\mfF := (\mcF, \nabla_\mcF)$ be a log integrable vector bundle on $\mfX_{/S}$.
If $\partial$ is a logarithmic derivation corresponding to a local section of $\mcT_{X^\mr{log}/S^\mr{log}}$,
then we shall denote by $\partial^{(p)}$ the $p$-th symbolic power of $\partial$ (i.e., ``$\partial \mapsto \partial^{(p)}$"  asserted  in  ~\cite{Og},  Proposition 1.2.1),
which is also a logarithmic derivation corresponding to a local section  of $\mcT_{X^\mr{log}/S^\mr{log}}$.
Then there exists uniquely an {\it $\mcO_X$-linear} morphism
\begin{equation}
{\psi}^{\nabla_\mcF} : \mcT_{X^\mr{log}/S^\mr{log}}^{\otimes p} \migi \mcE nd_{\mcO_X} (\mcF)
\end{equation}
determined  by assigning
\begin{equation}
\partial^{\otimes p} \mapsto \nabla_\mcF (\partial)^{\circ p} -  \nabla_\mcF (\partial^{(p)})
\end{equation}
for any local section $\partial \in \mcT_{X^\mr{log}/S^\mr{log}}$,  where $ \nabla_\mcF (\partial)^{\circ p}$ denotes the $p$-th iterate of the $f^{-1}(\mcO_S)$-linear endomorphism $\nabla_\mcF (\partial)$ of $\mcF$.
We shall refer to ${\psi}^{\nabla_\mcF}$ as   the {\it $p$-curvature} of $\nabla_\mcF$.

\vspace{5mm}
\subsection{} \label{09}

Next, we  recall  the monodromy of a logarithmic connection.
Let $(\mcF, \nabla_\mcF)$ be as above, and suppose that $r >0$.
For each $i =1, \cdots, r$, consider the composite
\begin{equation}
\mcF \stackrel{\nabla_\mcF}{\migi} \Omega_{X^\mr{log}/S^\mr{log}}\otimes \mcF  \migi
\sigma_{i*}(\sigma_i^{*}(\Omega_{X^\mr{log}/S^\mr{log}})\otimes \sigma_i^{*}(\mcF))
\isom \sigma_{i*}(\sigma_i^{*}(\mcF)),
\end{equation}
where the second arrow arises from 
the adjunction relation ``$\sigma_i^{*}(-) \dashv \sigma_{i*}(-)$"
(i.e., ``the functor $\sigma^{*}_i(-)$ is left adjoint to the functor $\sigma_{i*}(-)$"), and the third arrow arises   from the residue  isomorphism (\ref{triv}).
This composite corresponds (via the adjunction relation `$\sigma_i^{*}(-) \dashv \sigma_{i*}(-)$") to an  $\mcO_{ S}$-linear endomorphism  $\sigma_i^{*}(\mcF) \migi \sigma_i^{*}(\mcF)$, equivalently, 
 a global section
\begin{equation}
\mu_i^{ \nabla_\mcF} \in \Gamma ( S, \mcE nd_{\mcO_{ S}} (\sigma_i^{*}(\mcF))).
\end{equation}
\vspace{3mm}
\bde  \label{z010}\leavevmode\\
 \ \ \ 
 We shall refer to $\mu_i^{\nabla_\mcF}$ as the {\it monodromy} of $\nabla_\mcF$ at  $\sigma_i$.
\ede

\begin{rema}  \label{zz050} \leavevmode\\
 \ \ \ 
 Let $\mfL : = (\mcL, \nabla_\mcL)$ be a log integrable {\it line} bundle on $\mfX_{/S}$.
Then, $\mcE nd_{\mcO_S} (\sigma_i^*(\mcL)) \cong \mcO_S$, and hence, $\mu_i^{\nabla_\mcL}$ ($i = 1, \cdots, r$) may be thought of as an element of $\Gamma (S, \mcO_S)$.
In particular, it makes  sense to ask whether $\mu_i^{\nabla_\mcL}$ lies in $k$ ($\subseteq \Gamma (S, \mcO_S)$) or not.
 \end{rema}
\vspace{3mm}

\begin{rema}  \label{z08} \leavevmode\\
 \ \ \ If $\mcG$ is a vector bundle on $X^{(1)}_S$, then one may define (cf. \cite{Wak5}, \S\,3.3) canonically  an $S$-log connection
\begin{equation}
\nabla^\mr{can}_\mcG : F^*_{X/S}(\mcG) \migi \Omega_{X^\mr{log}/S^\mr{log}}\otimes F^*_{X/S}(\mcG) 
\end{equation}
on the pull-back $F^*_{X/S}(\mcG)$ of $\mcG$, which is uniquely determined by the condition that the sections of the subsheaf $F^{-1}_{X/S}(\mcG)$ ($\subseteq F^*_{X/S}(\mcG)$) are contained in $\mr{Ker}(\nabla^\mr{can}_\mcG)$.
We shall refer to $\nabla^\mr{can}_\mcG$ as the {\it canonical $S$-log connection} on $F^*_{X/S}(\mcG)$.
One verifies immediately  that
\begin{equation} \label{k09}
\mr{Im} (\nabla^\mr{can}_\mcG) \subseteq \Omega_{X/S}\otimes F^*_{X/S}(\mcG) \ \left(\subseteq \Omega_{X^\mr{log}/S^\mr{log}}\otimes F^*_{X/S}(\mcG) \right)
\end{equation}
(i.e., $\nabla_\mcG^\mr{can}$ arises from a {\it non-logarithmic} connection on $F^*_{X/S}(\mcG)$) and
\begin{equation} \label{k08}
\psi^{\nabla^\mr{can}_\mcG} = \mu_i^{\nabla_\mcG^\mr{can}} = 0
\end{equation}
for any $i$.

 
\end{rema}

\vspace{5mm}
\subsection{} \label{09090}

Let us write  $m := \mr{rk}(\mcF)$.
For each $i = 1, \cdots, r$, denote by $\phi_i^{ \nabla_\mcF} (t) \in \Gamma (S, \mcO_S) [t]$ the  characteristic polynomial of  $\mu_i^{\nabla_\mcF}$, i.e., 
\begin{equation} \label{k011}
\phi_i^{ \nabla_\mcF} (t) := \mr{det}(t \cdot \mr{id}_{\sigma_i^*(\mcF)}- \mu_i^{ \nabla_\mcF}) = \sum_{j =0}^m a^{\nabla_\mcF}_{i, j} \cdot t^j,
\end{equation}
where $a^{\nabla_\mcF}_{i, j} \in  \Gamma (S, \mcO_S) $ (satisfying that $a^{\nabla_\mcF}_{i, m} = 1$).

\vspace{3mm}
\bpr  \label{torsion} \leavevmode\\
 \ \ \ 
Suppose  further that ${ \psi}^{\nabla_\mcF} =0$.
Then, for any $j = 0, \cdots, m$, the element $a_{i, j}^{\nabla_\mcF}$ lies in $\mbF_p$.
  \epr
\begin{proof}
The condition ``${ \psi}^{\nabla_\mcF} =0$"  implies the equality 
\begin{equation}
(\mu_i^{\nabla_\mcF})^{p} -  \mu_i^{\nabla_\mcF} =0.
\end{equation}
Hence, (since $S$ is of characteristic $p$)  the following  sequence of equalities holds:
\begin{align}
\sum_{j =0}^m (a_{i, j}^{\nabla_\mcF})^p \cdot t^{j p} & = (\sum_{j =0}^m a_{i, j}^{\nabla_\mcF} \cdot t^j)^p \\
& =\mr{det} ((t \cdot \mr{id}_{\sigma^*_i(\mcF)}- \mu_i^{\nabla_\mcF})^p) \notag \\
& = \mr{det} (t^p \cdot  \mr{id}_{\sigma^*_i(\mcF)} - (\mu_i^{\nabla_\mcF})^p)\notag \\
& = \mr{det} (t^p \cdot \mr{id}_{\sigma^*_i(\mcF)} - \mu_i^{\nabla_\mcF}) \notag\\
& = \sum_{j =0}^m a_{i, j}^{\nabla_\mcF} \cdot t^{jp}.\notag
\end{align}
This yields the   equality 
\begin{equation}
(a_{i, j}^{\nabla_\mcF})^p =  a_{i, j}^{\nabla_\mcF}
\end{equation}
($j = 0, \cdots, m$), i.e., $a_{i, j}^{\nabla_\mcF} \in \mbF_p$.
This completes the proof of Proposition \ref{torsion}.
\end{proof}
\vspace{3mm}


\bde  \label{z012}\leavevmode\\
 \ \ \ 
Let $\mfF :=(\mcF, \nabla_\mcF)$ be as above.
\begin{itemize}
\item[(i)]
Let  $\tau_i := [\tau_{i, 1}, \cdots, \tau_{i, m}]$
 be  a multiset over $\Gamma (S, \mcO_S)$ with cardinality $m$.
We shall say that  {\it $\nabla_\mcF$ is of exponent $\tau_i$} at $\sigma_i$ if $\phi_i^{\nabla_\mcF} (t)$ may be described as
\begin{equation}
\phi_i^{\nabla_\mcF} (t) = \prod_{j =1}^m (t -\tau_{i, j}).
\end{equation} 
\item[(ii)]
Suppose that $r >0$ and we are given an $r$-tuple 
 $\vec{\tau} := (\tau_i)_{i=1}^r$ of multisets over $k$ with cardinality $n$ (i.e., an element of $(\mbN^k_{\sharp (-) = n})^{\times r}$).
 Then, we shall say that $\mfF$ is {\it of exponent $\vec{\tau}$} if $\nabla_\mcF$ is of exponent $\tau_i$ at $\sigma_i$ for any $i \in \{ 1, \cdots, r \}$.
\item[(iii)]
Suppose that $r =0$.
Then, we shall say, for convenience, that any log integrable vector bundle 
 on $\mfX_{/S}$ is {\it of exponent $\emptyset$}.
\end{itemize}
\ede
\vspace{3mm}

\begin{rema}  \label{z08s4} \leavevmode\\
 \ \ \ 
 Let $i \in \{ 1, \cdots, r \}$,  $\tau_i \in  \Gamma (S, \mcO_S)$, and let  $\mfL := (\mcL, \nabla_\mcL)$ be 
 a log integrable line bundle on $\mfX_{/S}$.
 Then, $\nabla_\mcL$ is of exponent $\tau_i$ at $\sigma_i$ if and only if $\mu_i^{\nabla_\mcL}=  \tau_i$.
  \end{rema}

\vspace{3mm}

The following two propositions follow immediately  from the various definitions involved.

\vspace{3mm}
\bpr  \label{torsion77} \leavevmode\\
 \ \ \ 
 Let  $\mfF = (\mcF, \nabla_\mcF)$ be   a log integrable vector bundle on $\mfX_{/S}$,
  and suppose that $\nabla_\mcF$ is  of exponent $\tau_i := [\tau_{i, 1}, \cdots, \tau_{i, m}]$ (where $\tau_{i, j} \in \Gamma (S, \mcO_S)$) at $\sigma_i$.
 Then, the $S$-log connection  $\mr{det}(\nabla_\mcF)$ on the line bundle $\mr{det}(\mcF)$ is of exponent $\sum_{j=1}^m \tau_{i, j}$ at $\sigma_i$.
  \epr

\vspace{3mm}
\bpr  \label{torsion7w7} \leavevmode\\
 \ \ \ 
 Let  $\mfF := (\mcF, \nabla_\mcF)$ be a log integrable vector bundle on $\mfX_{/S}$ of exponent $\vec{\tau} := (\tau_i)_{i=1}^r \in (\mbN^k_{\sharp (-) =n})^{\times r}$ and
$\mfL := (\mcL, \nabla_\mcL)$  a log integrable line bundle on $\mfX_{/S}$ of exponent $\vec{a} := (a_i)_{i =1}^r \in k^{\times r}$ ($= (\mbN^k_{\sharp (-) =1})^{\times r}$).
 Then, the log integrable vector  bundle 
 \begin{equation}
 \mfF \otimes \mfL := (\mcF \otimes \mcL, \nabla_\mcF \otimes \nabla_\mcL)
 \end{equation}
 is of exponent  $\vec{\tau}^{\, + \vec{a}}$ (cf. (\ref{z003})).
  \epr
\vspace{5mm}

\vspace{6mm}
\section{Determinant data} \label{z013} \vspace{3mm}

 In this section, we shall recall the notion of a determinant data (cf. Definition \ref{d7}, (i)), which was introduced  in the author's paper (cf. 
 ~\cite{Wak5}, Definition 4.9.1 (i)).
As we proved in ~\cite{Wak5} (cf. ~\cite{Wak5}, Theorem D, or (\ref{w5001}) displayed later),
one may realize, after fixing an $n$-determinant data $\mbU$, each $\mfs \mfl_n$-oper as a certain integrable vector bundle, i.e., a $(\mr{GL}_n, \mbU)$-oper.
Since we have assumed that the ground field $k$ is of positive characteristic, there exists (cf. Proposition \ref{z03055ww}) necessarily an $n$-determinant data with prescribed monodromy.

\vspace{5mm}
\subsection{} \label{z014}

Let $S$ and  $\mfX_{/S}$ be as before, and $n$ a positive  integer with $n \leq p$.
\vspace{3mm}
\bde  (cf. 
 ~\cite{Wak5}, Definition 4.9.1) \label{d7}\leavevmode\\
 \ \ \ 
 \vspace{-5mm}
\begin{itemize}
\item[(i)]
An {\it $n$-determinant data} for $\mfX_{/S}$ is a pair 
\begin{equation}
\mbU := (\BB, \nabla_0)
\end{equation}
 consisting of a line bundle $\BB$ on $X$ and an $S$-log connection $\nabla_0$ on 
 the line bundle $\mcT_{X^\mr{log}/S^\mr{log}}^{\otimes \frac{n (n-1)}{2}}\otimes \BB^{\otimes n}$.
\item[(ii)]
Let $\mbU := (\BB, \nabla_0)$ and $\mbU' := (\BB', \nabla'_0)$ be  $n$-determinant data
for $\mfX_{/S}$.
An {\it isomorphism of $n$-determinant data} from $\mbU$ to $\mbU'$ is an isomorphism
$\BB \isom \BB'$ of $\mcO_X$-modules such that 
the induced isomorphism
\begin{equation}
\mcT_{X^\mr{log}/S^\mr{log}}^{\otimes \frac{n (n-1)}{2}}\otimes \BB^{\otimes n} \isom \mcT_{X^\mr{log}/S^\mr{log}}^{\otimes \frac{n (n-1)}{2}}\otimes \BB'^{\otimes n}
\end{equation}
 is compatible with the respective $S$-log connections $\nabla_0$ and $\nabla'_0$.
\end{itemize}
  \ede

\begin{rema} \label{y002} \leavevmode\\
 \ \ \ Let $\mbU := (\BB, \nabla_0)$ be an $n$-determinant data for $\mfX_{/S}$ and  $s' : S' \migi S$  a morphism of $k$-schemes.
Then the base-change
\begin{equation} \label{PB}
s'^* (\mbU) := ((\mr{id}_X \times s')^*(\BB), (\mr{id}_X \times s')^* (\nabla_0))
\end{equation}
via $s'$ forms an $n$-determinant for $\mfX_{/S'}$ (cf. (\ref{BC44})).
\end{rema}
\vspace{3mm}

\bde \label{z015} \leavevmode\\
 \ \ \ 
We shall say that an
$n$-determinant data $\mbU := (\BB, \nabla_0)$ for $\mfX_{/S}$ is {\it dormant} if $\psi^{\nabla_0} =0$.
\ede
\bde \label{z015dd} \leavevmode\\
 \ \ \ 
Let $\vec{a}$ be an element of $k^{\times r}$. 
We shall say that an
$n$-determinant data $\mbU := (\BB, \nabla_0)$ for $\mfX_{/S}$ is {\it of exponent $\vec{a}$} if 
the log integrable line bundle $(\mcT_{X^\mr{log}/S^\mr{log}}^{\otimes \frac{n (n-1)}{2}}\otimes \BB^{\otimes n}, \nabla_0)$ is of exponent $\vec{a}$ (cf. Definition \ref{z012} (ii) and (iii)).
\ede

\vspace{3mm}
\bpr \label{z03055ww}\leavevmode\\
 \ \ \ 
Suppose that $r>0$, and let $\vec{a} := (a_i)_{i=1}^r \in \mbF_p^{\times r}$.
Then,   there exists a dormant  $n$-determinant data $\mbU := (\BB, \nabla_0)$  for $\mfX_{/S}$ of exponent $\vec{a}$.
\epr
\begin{proof}
Since $n \leq p$, one may choose a pair of nonnegative integers $(s, t)$ satisfying that 
$p \cdot s = n \cdot t + \frac{n (n-1)}{2}$.
Let us take $\BB' := \mcT_{X^\mr{log}/S^\mr{log}}^{\otimes t}$.
Then, 
\begin{equation}
\mcT_{X^\mr{log}/S^\mr{log}}^{\otimes \frac{n (n-1)}{2}}\otimes \BB'^{\otimes n} \isom (\mcT_{X^\mr{log}/S^\mr{log}}^{\otimes s})^{\otimes p} \isom F_{X/S}^*((\mr{id}_X \times F_S)^*(\mcT_{X^\mr{log}/S^\mr{log}}^{\otimes s})).
\end{equation}
Denote by $\nabla'_0$ the $S$-log connection on $\mcT_{X^\mr{log}/S^\mr{log}}^{\otimes \frac{n (n-1)}{2}}\otimes \BB'^{\otimes n}$ corresponding, via this composite isomorphism, to  the canonical $S$-log connection on $F_{X/S}^*((\mr{id}_X \times F_S)^*(\mcT_{X^\mr{log}/S^\mr{log}}^{\otimes s}))$ (cf. Remark \ref{z08}).
Then, the  pair $(\BB', \nabla'_0)$ forms a dormant $n$-determinant data for $\mfX_{/S}$ 
of exponent $(0, 0, \cdots, 0)$
(cf. (\ref{k08})).

Now,  for each $i =1, \cdots, r$,
we shall choose  a nonnegative integer $m_i$ such that $\overline{m_i \cdot n} = a_i$ in $\mbF_p$, where $\overline{m_i \cdot n}$ denotes the image of $m_i \cdot n \in \mbZ$ via the quotient $\mbZ \migisurj \mbF_p$.
In particular, $\sum_{i=1}^r (m_i \cdot n) \sigma_i$ is an effective relative divisor on $X$ relative to $S$.
By passing to  the isomorphism
\begin{equation}
(\mcT_{X^\mr{log}/S^\mr{log}}^{\otimes \frac{n (n-1)}{2}}\otimes \BB'^{\otimes n}) (-\sum_{i=1}^r (m_i \cdot n) \sigma_i) \isom \mcT_{X^\mr{log}/S^\mr{log}}^{\otimes \frac{n (n-1)}{2}}\otimes \BB' (-\sum_{i=1}^r m_i \sigma_i)^{\otimes n}, 
\end{equation}
we obtain an $S$-log connection  $\nabla_0$ on $\mcT_{X^\mr{log}/S^\mr{log}}^{\otimes \frac{n (n-1)}{2}}\otimes \BB' (-\sum_{i=1}^r m_i \sigma_i)^{\otimes n}$  (with vanishing $p$-curvature)  corresponding to 
the restriction of $\nabla'_0$ to $(\mcT_{X^\mr{log}/S^\mr{log}}^{\otimes \frac{n (n-1)}{2}}\otimes \BB'^{\otimes n})( -\sum_{i=1}^r (m_i \cdot n) \sigma_i)$.
This $S$-log connection is, by construction, of exponent $\overline{m_i \cdot n} =  a_i$ at $\sigma_i$.
Thus, we obtain a dormant  $n$-determinant data
\begin{equation}
(\BB := \BB' (-\sum_{i=1}^r m_i \sigma_i), \nabla_0)
\end{equation}
which  satisfies the required conditions, as desired.
\end{proof}

\vspace{5mm}
\subsection{} \label{z016}

 Let $\mbU := (\BB, \nabla_0)$ be an $n$-determinant data for $\mfX_{/S}$ and 
$\mfL :=(\mcL, \nabla_\mcL)$ be a log integrable line bundle on $\mfX_{/S}$.
We shall consider the pair
\begin{equation} \label{w00987}
\mbU \otimes \mfL := (\BB \otimes \mcL, \nabla_0 \otimes \nabla_\mcL^{\otimes n}),
\end{equation}
where we regard $\nabla_0 \otimes \nabla_\mcL^{\otimes n}$ (cf. \S\,\ref{z06}) as an $S$-log connection
on 
\begin{equation}
\mcT_{X^\mr{log}/S^\mr{log}}^{\otimes \frac{n (n-1)}{2}}\otimes (\BB\otimes \mcL)^{\otimes n}   \ \left(\cong (\mcT_{X^\mr{log}/S^\mr{log}}^{\otimes \frac{n (n-1)}{2}}\otimes \BB^{\otimes n}) \otimes \mcL^{\otimes n}\right).
\end{equation}
One verifies that $\mbU \otimes \mfL$ forms an $n$-determinant data for $\mfX_{/S}$.
If, moreover, $\mbU$ is dormant and $\psi^{\nabla_\mcL} =0$, then  $\mbU \otimes \mfL$ turns out to be  dormant.

\vspace{5mm}
\subsection{} \label{z018}

For each line bundle $\BB$ on $X$, one may construct, in a canonical manner,  
a dormant  $p$-determinant data
whose underlying line bundle coincides with $\BB$.

Indeed, let us consider the natural  composite isomorphism
\begin{equation}
F^*_{X/S}((\mr{id}_X \times F_S)^*(\mcT_{X^\mr{log}/S^\mr{log}}^{\otimes \frac{p-1}{2}}\otimes\BB)) \isom  (\mcT_{X^\mr{log}/S^\mr{log}}^{\otimes \frac{p-1}{2}}\otimes\BB)^{\otimes p}  \isom \mcT_{X^\mr{log}/S^\mr{log}}^{\otimes \frac{p(p-1)}{2}}\otimes\BB^{\otimes p}.
\end{equation}
By passing to this composite, 
we obtain an $S$-log connection $\nabla_{0, \BB}^\mr{can}$  on $\mcT_{X^\mr{log}/S^\mr{log}}^{\otimes \frac{p(p-1)}{2}}\otimes\BB^{\otimes p}$ 
corresponding to 
the canonical $S$-log connection 
on $F^*_{X/S}((\mr{id}_X \times F_S)^*(\mcT_{X^\mr{log}/S^\mr{log}}^{\otimes \frac{p-1}{2}}\otimes\BB))$ (cf. Remark \ref{z08}).
Thus, the pair
\begin{equation} \label{y09}
\mbU_\BB^\mr{can} := (\BB, \nabla_{0, \BB}^\mr{can})
\end{equation} 
forms a dormant $p$-determinant data for $\mfX_{/S}$ satisfying that $\mu_i^{\nabla_{0, \BB}^\mr{can}} = 0$ for any $i$ (cf. (\ref{k08})).
\vspace{5mm}
\subsection{} \label{z019}

Let $\mbU := (\BB, \nabla_0)$ be an $n$-determinant data  for $\mfX_{/S}$, and write  
\begin{equation} \label{y01}
\BB^{\triangledown, n} := \Omega_{X^\mr{log}/S^\mr{log}}^{\otimes (n-1)}\otimes \BB^\vee.
\end{equation}
 Then, we have a  natural composite  isomorphism
\begin{equation}
\mcT_{X^\mr{log}/S^\mr{log}}^{\otimes \frac{n(n-1)}{2}}\otimes (\BB^{\triangledown, n})^{\otimes n} \isom \mcT_{X^\mr{log}/S^\mr{log}}^{\otimes \frac{n(n-1)}{2}}\otimes (\Omega_{X^\mr{log}/S^\mr{log}}^{\otimes (n-1)}\otimes \BB^\vee)^{\otimes n} \isom (\mcT_{X^\mr{log}/S^\mr{log}}^{\otimes \frac{n(n-1)}{2}}\otimes \BB^{\otimes n})^\vee.
\end{equation}
The $S$-log connection $\nabla_0^\vee$ on $(\mcT_{X^\mr{log}/S^\mr{log}}^{\otimes \frac{n(n-1)}{2}}\otimes \BB^{\otimes n})^\vee$ carries,  by means of  this composite isomorphism, an $S$-log connection on $\mcT_{X^\mr{log}/S^\mr{log}}^{\otimes \frac{n(n-1)}{2}}\otimes (\BB^{\triangledown, n})^{\otimes n}$; we shall denote this connection by $\nabla^\triangledown_0$.
Thus, we obtain an $n$-determinant data
\begin{equation} \label{y02}
\mbU^\triangledown := (\BB^{\triangledown, n}, \nabla_0^\triangledown)
\end{equation}
 for $\mfX_{/S}$, which is referred to as
 the {\it dual $n$-determinant data} of $\mbU$.
One verifies that there exists a natural isomorphism
\begin{equation}
(\mbU^\triangledown)^\triangledown \isom \mbU
\end{equation}
of $n$-determinant data. 

\vspace{5mm}
\subsection{} \label{z029}

Let $\mbU  := (\BB, \nabla_0)$ be as above, and
write 
\begin{equation} \label{y03}
\BB^{\, \triangleright, n} := \mcT_{X^\mr{log}/S^\mr{log}}^{\otimes n} \otimes \BB.
\end{equation}
Consider the canonical composite isomorphism
\begin{align}
& \ \ \ \  \   (\mcT_{X^\mr{log}/S^\mr{log}}^{\otimes \frac{p (p-1)}{2}}\otimes \BB^{\otimes p})\otimes (\mcT_{X^\mr{log}/S^\mr{log}}^{\otimes \frac{n (n-1)}{2}}\otimes \BB^{\otimes n})^\vee \\
& \isom \mcT_{X^\mr{log}/S^\mr{log}}^{\otimes \frac{1}{2}\cdot (p^2 -p-n^2 +n)}  \otimes \BB^{\otimes (p-n)} \notag \\
& \isom \mcT_{X^\mr{log}/S^\mr{log}}^{\otimes \frac{(p-n) (p-n-1)}{2}} \otimes (\BB^{\triangleright,n})^{\otimes (p-n)}. \notag
\end{align}
The product $\nabla_{0, \BB}^\mr{can}\otimes \nabla_0^\vee$ of the $S$-log connections  carries, by means of this composite isomorphism, an $S$-log connection $\nabla_0^\triangleright$ on $\mcT_{X^\mr{log}/S^\mr{log}}^{\otimes \frac{(p-n) (p-n-1)}{2}} \otimes (\BB^{\triangleright, n})^{\otimes (p-n)}$.
Thus, we obtain a $(p-n)$-determinant data
\begin{equation} \label{y04}
\mbU^\triangleright := (\BB^{\triangleright, n}, \nabla_0^\triangleright)
\end{equation}
for $\mfX_{/S}$.

Moreover, we shall write
\begin{equation} \label{y2345}
\mbU^\bigstar := (\mbU^\triangleright)^\triangledown
\end{equation}
and  refer to it as the {\it $\bigstar$-dual $(p-n)$-determinant data} of $\mbU$.
If $\BB^\bigstar$ denotes 
the underlying line bundle  of $\mbU^\bigstar$, i.e.,
\begin{equation}
\BB^\bigstar := (\BB^{\triangleright, n})^{\triangledown, p-n},
\end{equation}
 then there exist  natural isomorphisms
\begin{equation} \label{y0012}
\BB^\bigstar \isom \Omega^{\otimes (p-1)}_{X^\mr{log}/S^\mr{log}} \otimes \BB^\vee \ \ \text{and} \ \ \BB^\bigstar \isom \BB^{\triangledown, p}.
\end{equation}

The following three propositions follow immediately from the various definitions involved.

\vspace{3mm}
\bpr \label{z030}\leavevmode\\
 \ \ \ 
There exists a canonical isomorphism 
\begin{equation}
(\mbU^\bigstar)^\bigstar  \isom \mbU
\end{equation}
of $n$-determinant data.
\epr

\vspace{3mm}
\bpr \label{z0030}\leavevmode\\
 \ \ \ 
Let $\mfL : = (\mcL, \nabla_\mcL)$ be a log integrable line bundle on $\mfX_{/S}$, and write
$\mfL^\vee := (\mcL^\vee, \nabla_\mcL^\vee)$.
Then,  there exists a canonical isomorphism 
\begin{equation}
(\mbU \otimes \mfL)^\bigstar  \isom \mbU^\bigstar \otimes \mfL^\vee
\end{equation}
of $(p-n)$-determinant data.
\epr

\vspace{3mm}
\bpr \label{z03055}\leavevmode\\
 \ \ \ 
Let $\BB$ be a line bundle on $X$.
Then, there exists a canonical isomorphism
\begin{equation}
\left(\mbU^\mr{can}_{\BB^{\triangledown, p}}  \isom \right) \ \mbU^\mr{can}_{\BB^{\bigstar}} \isom (\mbU^\mr{can}_\BB)^\triangledown
\end{equation}
of $p$-determinant data for $\mfX_{/S}$.
\epr

\vspace{6mm}
\section{Opers on pointed stable curves}  \label{z040}\vspace{3mm}

In this section, we recall the definition of a (dormant) $\mr{GL}_n$-oper  (where $\mr{GL}_n$ denotes the general linear group over $k$ of rank $n$) and consider the canonical construction 
 of  a dormant $\mr{GL}_p$-oper  by means of a line bundle (cf. Proposition \ref{778d3}).

Let $S$, $\mfX_{/S}$, and  $n$ be as before.
\vspace{5mm}
\subsection{}
First, recall the definition of a $\mr{GL}_n$-oper, as follows.
\vspace{3mm}
\bde (cf. ~\cite{Wak5}, Definition 4.2.1)  \label{z041}\leavevmode\\
\vspace{-5mm}
\begin{itemize}
 \item[(i)]
 A {\it $\mr{GL}_n$-oper} on  $\mfX_{/ S}$ is   a collection of  data 
 \begin{equation} 
 \label{GL1}
 \mfF^\heartsuit := (\mcF, \nabla_\mcF, \{ \mcF^j \}_{j=0}^n),\end{equation}
where
\begin{itemize}
\item[$\bullet$]
$\mcF$ is a vector bundle  on $X $ of rank $n$;
\item[$\bullet$]
$\nabla_\mcF$ is an $S$-log connection
 $\mcF \migi \Omega_{X^\mr{log}/S^\mr{log}} \otimes\mcF$ on $\mcF$;
\item[$\bullet$]
$\{ \mcF^j \}_{j=0}^n$ is   a decreasing filtration
\begin{equation}
0 = \mcF^n \subseteq \mcF^{n-1} \subseteq \dotsm \subseteq  \mcF^0= \mcF
\end{equation}
 on $\mcF$ by vector bundles on $X$,
\end{itemize}
 satisfying the following three conditions:
\vspace{1mm}
\begin{itemize}
\item[(1)]
The subquotients $\mcF^j / \mcF^{j+1}$ ($ 0\leq j\leq n-1$)  are line bundles;
\vspace{1mm}
\item[(2)]
$\nabla_\mcF (\mcF^j) \subseteq  \Omega_{X^\mr{log}/S^\mr{log}} \otimes \mcF^{j-1}$ ($1 \leq j \leq n-1$);
\vspace{1mm}
\item[(3)]
The {\it $\mcO_X$-linear} morphism
\begin{equation} \label{GL2}
\mfk \mfs^j_{\mfF^\heartsuit} :  \mcF^j/\mcF^{j+1} \stackrel{}{\migi}   \Omega_{X^\mr{log}/S^\mr{log}} \otimes (\mcF^{j-1}/\mcF^j)\end{equation}
defined by assigning $\overline{a} \mapsto \overline{\nabla_\mcF (a)}$ for any local section $a \in \mcF^j$ (where $\overline{(-)}$'s denote the images in the respective quotients), which is well-defined by virtue of the condition (2),  is an isomorphism.
\end{itemize}


\item[(ii)]
Let $\mfF^\heartsuit := (\mcF, \nabla_\mcF, \{ \mcF^j \}_{j=0}^n)$, $\mfG^\heartsuit := (\mcG, \nabla_\mcG, \{ \mcG^j \}_{j=0}^n)$ be $\mr{GL}_n$-opers on $\mfX_{/S}$.
An {\it isomorphism of $\mr{GL}_n$-opers} from $\mfF^\heartsuit$ to $\mfG^\heartsuit$ is an isomorphism
$(\mcF, \nabla_\mcF) \isom (\mcG, \nabla_\mcG)$ of log integrable vector bundles (cf. \S\,\ref{z06})  that is compatible with the respective filtrations $\{ \mcF^j\}_{j=0}^n$ and  $\{ \mcG^j \}_{j=0}^n$.
\end{itemize}  
  \ede

\begin{rema} \label{z042} \leavevmode\\
 \ \ \ Let $\mfF^\heartsuit := (\mcF, \nabla_\mcF, \{ \mcF^j \}_{j =0}^n)$ be a $\mr{GL}_n$-oper on $\mfX_{/S}$ and fix $j  \in \{ 0, \cdots, n-1\}$.
By composing  the  isomorphisms 
\begin{equation}
 (\mcF^{l}/\mcF^{l+1}) \otimes  \mcT_{X^\mr{log}/S^\mr{log}}^{\otimes (l-j)} \isom  (\mcF^{l+1}/\mcF^{l+2}) \otimes  \mcT_{X^\mr{log}/S^\mr{log}}^{\otimes (l-j +1)}
 \end{equation}
arising from the isomorphisms 
$\mfk \mfs_{\mfF^\heartsuit}^l$ ($j+1\leq l \leq  n-1$) (cf.  (\ref{GL2})),  we obtain a canonical composite isomorphism 
\begin{equation} \label{d13}  \mcF^j/\mcF^{j+1} \isom (\mcF^{j+1}/\mcF^{j+2}) \otimes  \mcT_{X^\mr{log}/S^\mr{log}} \isom \cdots \isom \mcF^{n-1} \otimes  \mcT_{X^\mr{log}/S^\mr{log}}^{\otimes n-1-j} \end{equation}
between line bundles on $X$.
Moreover, by using  these isomorphisms for all $j$, we obtain a  composite isomorphism 
\begin{align} \label{GL30}
 \mfd \mfe \mft_{\mfF^\heartsuit} : 
 \mr{det}(\mcF) & \isom \bigotimes_{j=0}^{n-1}(\mcF^j/\mcF^{j+1}) \\
 &  \isom \big(\bigotimes_{j=0}^{n-1}\mcT_{X^\mr{log}/S^\mr{log}}^{\otimes (n-1-j)} \big) \otimes  (\mcF^{n-1})^{\otimes n}\notag \\
 &  \isom \mcT_{X^\mr{log}/S^\mr{log}}^{\otimes \frac{n(n-1)}{2}} \otimes  (\mcF^{n-1})^{\otimes n}. 
\notag \end{align}
In particular, the following equalities hold:
\begin{align} \label{y13}
\mr{deg} (\mcF) & = \mr{deg} (\mcT_{X^\mr{log}/S^\mr{log}}^{\otimes \frac{n(n-1)}{2}} \otimes  (\mcF^{n-1})^{\otimes n})  = n (n-1) (1-g) + n \cdot \mr{deg} (\mcF^{n-1}).
\end{align}

\end{rema}
\vspace{5mm}

The following proposition will be used  in the proof of Proposition \ref{z07311}.

\vspace{3mm}
\bpr \label{d3} \leavevmode\\
 \ \ \ 
Let $\mfF^\heartsuit := (\mcF, \nabla_\mcF, \{ \mcF^j \}_{j=0}^n)$ be  a $\mr{GL}_n$-oper on $\mfX_{/S}$ and $(\mcV, \nabla_\mcV)$  a log integrable vector bundle on $\mfX_{/S}$. 
Suppose that we are given two morphisms $\phi_1$, $\phi_2 : \mcF \migi \mcV$ both of which are compatible with the respective $S$-log connections $\nabla_\mcF$, $\nabla_\mcV$, and satisfying that $\phi_1 |_{\mcF^{n-1}} = \phi_2 |_{\mcF^{n-1}}$.
Then, the equality $\phi_1 = \phi_2$ holds. 
\epr
\begin{proof}
Suppose that $\phi_1 |_{\mcF^{j}} = \phi_2 |_{\mcF^{j}}$ for some $j \in \{ 1, \cdots , n-1 \}$.
By the  definition of  a $\mr{GL}_n$-oper,  
 $\mcF^{j-1}$ may be  generated  by $\mcF^{j}$ and $\nabla_\mcF (\mcF^{j})$.
 Hence, since both $\phi_1$ and $\phi_2$ are compatible with the $S$-log connections $\nabla_\mcF$ and $\nabla_\mcV$, 
 the equality $\phi_1 |_{\mcF^{j}} = \phi_2 |_{\mcF^{j}}$ implies the equality 
 $\phi_1 |_{\mcF^{j-1}} = \phi_2 |_{\mcF^{j-1}}$.
 Thus, the assertion follows from descending induction on $j$.
\end{proof}

\vspace{3mm}
\bde \label{z044}\leavevmode\\
 \ \ \ 
We shall say that a $\mr{GL}_n$-oper $\mfF^\heartsuit := (\mcF, \nabla_\mcF, \{ \mcF^j \}_{j=0}^n)$ on $\mfX_{/S}$ is {\it dormant} if ${\psi}^{\nabla_\mcF}=0$.
\ede

\vspace{5mm}
\subsection{} \label{z045}

Let $\mfF := (\mcF, \nabla_\mcF)$ be  a log integrable vector bundle on $\mfX_{/S}$ of rank $n$.
We shall write 
\begin{equation} 
\label{1070}
\mcA_{\mr{Ker}(\nabla_\mcF)} := F^*_{X/S}(F_{X/S*}(\mr{Ker}(\nabla_\mcF))).\end{equation}
Note that although $\mr{Ker}(\nabla_\mcF)$ is not an $\mcO_X$-module, 
one may equip, in a natural manner,  its direct image $F_{X/S*}(\mr{Ker}(\nabla_\mcF))$ via $F_{X/S}$ with a structure of $\mcO_{X^{(1)}_S}$-module.
The  $\mcO_{X^{(1)}_S}$-linear   inclusion $F_{X/S*}(\mr{Ker}(\nabla_\mcF)) \migiincl F_{X/S*}(\mcF)$ corresponds,
via  the adjunction relation 
 ``$F^*_{X/S} (-) \dashv F_{X/S*}(-)$" ,
to an $\mcO_X$-linear  morphism 
\begin{equation} \label{inj77}
 \nu^{\nabla_\mcF} : \mcA_{\mr{Ker}(\nabla_\mcF)}  \migi \mcF.\end{equation}
If we consider the canonical $S$-log connection $\nabla_{F_{X/S*}(\mr{Ker}(\nabla_\mcF))}^{\mr{can}}$ on $\mcA_{\mr{Ker}(\nabla_\mcF)}$ (cf. Remark \ref{z08}), then the morphism $\nu^{\nabla_\mcF}$ is compatible with the respective connections $\nabla_{F_{X/S*}(\mr{Ker}(\nabla_\mcF))}^{\mr{can}}$ and $\nabla_\mcF$.
The morphism $\nu^{\nabla_\mcF}$ fits into the short exact sequence
\begin{equation}
0 \migi  \mcA_{\mr{Ker}(\nabla_\mcF)}  \stackrel{ \nu^{\nabla_\mcF}}{\migi} \mcF \stackrel{}{\migi} \Lambda_{\mr{sing}}\oplus \bigoplus_{i=1}^r \Lambda_i \migi 0
\end{equation}
of $\mcO_X$-modules (cf. ~\cite{Wak5},  the short exact sequence (872)), where $\Lambda_{\mr{sing}}$ and $\Lambda_i$ ($i =1, \cdots, r$) are  $\mcO_X$-modules supported on the nonsmooth   locus of  the semistable curve $X$ (over $S$) and the locus $\mr{Im}(\sigma_i) \subseteq X$ respectively.

Let us  recall the following proposition.

\vspace{3mm}
\bpr  \label{z046}\leavevmode\\
 \ \ \ 
 In the above notation, suppose further that $\psi^{\nabla_\mcF} = 0$  and $S = \mr{Spec}(k')$ for some algebraically closed field $k'$ over $k$. Then,  for each $i$ ($= 1, \cdots, r$),
 there exists uniquely a multiset $[\tau_{i, 1}, \cdots, \tau_{i, n}]$ over $\mbZ$
  with cardinality $n$ satisfying the following three  conditions:
 \begin{itemize}
 \item[(1)]
 $0 \leq \tau_{i, j} <p$ for any $j = 1, \cdots, n$;
  \item[(2)]
 $\Lambda_i$ decomposes into the direct sum 
\begin{equation}
\Lambda_i \isom \bigoplus_{j =1}^n \mcO_X (\tau_{i, j} \sigma_i)/\mcO_X;
\end{equation}
 \item[(3)]
$\nabla_\mcF$ is of exponent $[\overline{\tau}_{i, 1}, \cdots, \overline{\tau}_{i, n}]$ at $\sigma_i$  (cf. Definition \ref{z012} (i)), where 
 each $\overline{\tau}_{i, j}$  denotes the image of $\tau_{i, j}$ via the quotient  $\mbZ \migi \mbF_p$ ($\subseteq k$).
 \end{itemize}


\epr
\begin{proof}
The assertion follows from ~\cite{O2}, Corollary 2.10, or the discussion in ~\cite{Wak5} following Lemma 8.3.2.
\end{proof}

\vspace{5mm}
\subsection{} \label{z048}

Next, we consider $\mr{GL}_n$-opers with  prescribed determinant.
Let $\mbU := (\BB, \nabla_0)$ be an $n$-determinant data (cf. Definition \ref{d7}, (i)) for $\mfX_{/S}$.

\vspace{3mm}
\bde  (cf. ~\cite{Wak5}, Definition 4.9.4)\label{d1}\leavevmode\\
 \vspace{-5mm}
\begin{itemize}
\item[(i)]
A {\it $(\mr{GL}_n, \mbU)$-oper} on $\mfX_{/S}$ is 
a collection of data
\begin{equation}
\mfF^\diamondsuit := (\mcF, \nabla_\mcF, \{ \mcF^j \}_{j=0}^n, \eta_{\mfF^\diamondsuit}),
\end{equation}
where 
\begin{itemize}
\item[$\bullet$]
the collection of data 
\begin{equation}
\mfF^{\diamondsuit \heartsuit} := (\mcF, \nabla_\mcF, \{ \mcF^j \}_{j=0}^n)
\end{equation}
 is 
a $\mr{GL}_n$-oper  on $\mfX_{/S}$ (which is referred to  as the {\it underlying $\mr{GL}_n$-oper of $\mfF^\diamondsuit$});
\item[$\bullet$]
 $\eta_{\mfF^\diamondsuit}$ is an isomorphism $\BB \isom \mcF^{n-1}$ of $\mcO_X$-modules such that the composite isomorphism
\begin{equation} \label{y18}
 \mcT_{X^\mr{log}/S^\mr{log}}^{\otimes \frac{n (n-1)}{2}}\otimes  \mcB^{\otimes n}
\stackrel{\mr{id} \otimes \eta_{\mfF^\diamondsuit}^{\otimes n}}{\migi}
\mcT_{X^\mr{log}/S^\mr{log}}^{\otimes \frac{n (n-1)}{2}}\otimes (\mcF^{n-1})^{\otimes n} 
\stackrel{\mfd \mfe \mft^{-1}_{\mfF^{\diamondsuit \heartsuit} }}{\migi} 
\mr{det}(\mcF) 
\end{equation}
(cf. (\ref{GL30}) for the definition of the isomorphism $\mfd \mfe \mft_{(-)}$)
  is compatible with $\nabla_0$ and  $\mr{det}(\nabla_\mcF)$ (cf. \S\,\ref{z06}).
\end{itemize}

\item[(ii)]
Let $\mfF^\diamondsuit := (\mcF, \nabla_\mcF, \{ \mcF^j \}_{j=0}^n, \eta_{\mfF^\diamondsuit})$ and $\mfG^\diamondsuit := (\mcG, \nabla_\mcG, \{ \mcG^j \}_{j=0}^n, \eta_{\mfG^\diamondsuit})$ be 
$(\mr{GL}_n, \mbU)$-opers on $\mfX_{/S}$. 
An {\it isomorphism of $(\mr{GL}_n, \mbU)$-opers} from $\mfF^\diamondsuit$ to $\mfG^\diamondsuit$ is 
an isomorphism 
$\alpha : \mfF^{\diamondsuit \heartsuit} \isom \mfG^{\diamondsuit \heartsuit}$ between the respective underlying  
$\mr{GL}_n$-opers (cf. Definition \ref{z041}, (ii)) whose restriction  $\alpha |_{\mcF^{n-1}} : \mcF^{n-1} \isom \mcG^{n-1}$ satisfies the equality $\alpha |_{\mcF^{n-1}}  \circ \eta_{\mfF^\diamondsuit} = \eta_{\mfG^\diamondsuit}$.
\end{itemize}

  \ede
\vspace{3mm}

\bde \label{z0512}\leavevmode\\
 \ \ \ 
 Let  $\vec{\tau} := (\tau_i)_{i =1}^r$ be  an element of $(2^k_{\sharp (-) =n})^{\times r}$ (where we take $\vec{\tau} := \emptyset$ if $r =0$).
  Then, we shall say that a $(\mr{GL}_n, \mbU)$-oper $\mfF^\diamondsuit := (\mcF, \nabla_\mcF, \{ \mcF^j \}_{j=0}^n, \eta_{\mfF^\diamondsuit})$ on $\mfX_{/S}$ is {\it of exponent $\vec{\tau}$} if 
  the underlying  log integrable vector bundle $(\mcF, \nabla_\mcF)$ is of exponent $\vec{\tau}$ (cf. Definition \ref{z012} (ii) and (iii)).
\ede

\begin{rema} \label{z049} \leavevmode\\
 \ \ \ 
 Notice that the definition of a $(\mr{GL}_n, \mbU)$-oper in Definition \ref{d1}  differs slightly  from the definition of a $(\mr{GL}_n, 1, \mbU)$-oper proposed in ~\cite{Wak5}, Definition 4.9.4.
\end{rema}
\vspace{5mm}

\begin{rema} \label{z049ii} \leavevmode\\
 \ \ \ 
 Let $s' : S' \migi S$ be a morphism of $k$-schemes and $\mfF^\diamondsuit  := (\mcF, \nabla_\mcF, \{ \mcF^j \}_{j=0}^n, \eta_{\mfF^\diamondsuit})$ a $(\mr{GL}_n, \mbU)$-oper on $\mfX_{/S}$.
 Then, the collection of data
 \begin{equation} \label{w0090909}
 s'^* (\mfF^\diamondsuit)
 \end{equation}
 obtained from pulling-back the collection of data $\mfF^\diamondsuit$ via $\mr{id}_X \times s' : X \times_S S' \migi X$ forms a $(\mr{GL}_n, s'^*(\mbU))$-oper (cf. (\ref{PB})) on $\mfX_{/S'}$ (cf. (\ref{BC44})).
 \end{rema}
\vspace{5mm}

\bde \label{z051}\leavevmode\\
 \ \ \ 
We shall say that a $(\mr{GL}_n, \mbU)$-oper $\mfF^\diamondsuit  := (\mcF, \nabla_\mcF, \{ \mcF^j \}_{j=0}^n, \eta_{\mfF^\diamondsuit})$
  is {\it dormant} if $\psi^{\nabla_\mcF} =0$.
\ede

\begin{rema} \label{z04eer9} \leavevmode\\
 \ \ \ 
 If there exists  a dormant  $(\mr{GL}_n, \mbU)$-oper  $\mfF^\diamondsuit  := (\mcF, \nabla_\mcF, \{ \mcF^j \}_{j=0}^n, \eta_{\mfF^\diamondsuit})$ on $\mfX_{/S}$, then $\mbU$ is necessarily dormant (by  ~\cite{Wak5}, Proposition 3.2.2).
 Indeed, if $\mr{Tr}$ denotes  the trace map $\mr{Tr} : \mcE nd (\mcF) \migi \mcO_X$, then  we have
$\psi^{\nabla_0} =  \psi^{\mr{det} (\nabla_\mcF)} = \mr{Tr} \circ \psi^{\nabla_\mcF} =0$ (cf.
  ~\cite{JP}, Proposition 2.1.2 (iii)).
 \end{rema}
\vspace{5mm}

\begin{rema} \label{z04eer9g} \leavevmode\\
 \ \ \ 
 Suppose that $n =1$.
 Let us consider the $1$-step filtration $\{ \mcB^j \}_{j =0}^1$ on $\mcB$ given by 
 $\mcB^0 := \mcB$ and $\mcB^1 :=0$.
 Then,  the collection of data 
 \begin{equation} \label{sss}
\mfB^\diamondsuit :=  (\mcB, \nabla_0, \{ \mcB^j \}_{j=0}^1, \mr{id}_{\mcB})
 \end{equation}
 forms a unique (up to isomorphism) $(\mr{GL}_1, \mbU)$-oper on $\mfX_{/S}$.
If $\mbU$ is dormant, then $\mfB^\diamondsuit$ is tautologically dormant.
 \end{rema}
\vspace{5mm}

\vspace{3mm}
\bpr \label{z047} \leavevmode\\
 \ \ \ 
Suppose that $r >0$.
Let $\mbU := (\BB, \nabla_0)$ be an $n$-determinant data  for $\mfX_{/S}$ and 
$\mfF^\diamondsuit  := (\mcF, \nabla_\mcF, \{ \mcF^j \}_{j=0}^n, \eta_{\mfF^\diamondsuit})$
a dormant  $(\mr{GL}_n,  \mbU)$-oper on $\mfX_{/S}$ (hence $\mbU$ is dormant, as we explained in  Remark \ref{z04eer9}).
Then,  there exists an element $\vec{\tau} := (\tau_i)_{i=1}^r$ of $(2^{\mbF_p}_{\sharp (-) =n})^{\times r}$ such that $\mfF^\diamondsuit$ is of exponent $\vec{\tau}$.

\epr
\begin{proof}
It follows from Proposition \ref{torsion} that for each $i = 1, \cdots r$, the characteristic polynomial $\phi^{\nabla_\mcF}_i (t)$ of $\mu_i^{\nabla_\mcF}$  lies in $\mbF_p [t]$.
Thus, we may assume, after possibly restricting $X$ to a geometric fiber of $f : X \migi S$,  that $S = k'$ for some algebraically closed field $k'$ over $k$.
Let $[\tau_{i, 1}, \cdots, \tau_{i, n}]$ ($i =1, \cdots, r$) be the multiset asserted in Proposition \ref{z046} of the case where the log integrable vector bundle under consideration is taken to be $(\mcF, \nabla_\mcF)$. (In particular, $\nabla_\mcF$ is of exponent $\tau_i := [\overline{\tau}_{i, 1}, \cdots, \overline{\tau}_{i, n}]$ at $\sigma_i$.)
But, according to  ~\cite{Wak5}, Proposition 8.3.3, the integers $\tau_{i, 1}, \cdots, \tau_{i, n}$ are mutually distinct.
This implies that $\tau_i$ is consequently  a subset of $\mbF_p$,  
and hence,  completes the proof of Proposition \ref{z047}.
\end{proof}

\vspace{5mm}
\subsection{} \label{z0301}

Let $\mfF^\diamondsuit := (\mcF, \nabla_\mcF, \{ \mcF^j \}_{j=0}^n, \eta_{\mfF^\diamondsuit})$ be a $(\mr{GL}_n, \mbU)$-oper on $\mfX_{/S}$ of exponent $\vec{\tau} \in (\mbN^{k}_{\sharp (-) = n})^{\times r}$,  $\mfL : = (\mcL, \nabla_\mcL)$  a log integrable line bundle on $\mfX_{/S}$ of exponent   $\vec{a} := (a_i)_{i=1}^r \in k^{\times r}$.
Then, 
 one verifies from Proposition \ref{torsion7w7}
 that  the  collection of data
\begin{equation} \label{e00908}
\mfF^\diamondsuit_{\otimes \mfL} := (\mcF \otimes \mcL, \nabla_\mcF \otimes \nabla_\mcL, \{ \mcF^j \otimes \mcL\}_{j=0}^n, \eta_{\mfF^{\diamondsuit}} \otimes \mr{id}_\mcL)
\end{equation}
forms a $(\mr{GL}_n, \mbU \otimes \mfL)$-oper (cf. (\ref{w00987})) on $\mfX_{/S}$ of exponent $\vec{\tau}^{\,+ \vec{a}}$ (cf. (\ref{z003})).
\vspace{5mm}

\subsection{} \label{z0201}
Let $\mfF^\diamondsuit$ be as above. 
If we write $(\mcF^\vee)^j$ for the kernel of the dual $\mcF^\vee \migisurj (\mcF^{n-j})^\vee$ of the inclusion $\mcF^{n-j} \migiincl \mcF$,
then the collection of data
\begin{equation}
(\mcF^\vee, \nabla^\vee_\mcF, \{ (\mcF^\vee)^j \}_{j=0}^n)
\end{equation}
(cf.  \S\,\ref{z06} for the definition of $\nabla^\vee_\mcF$) forms a $\mr{GL}_n$-oper on $\mfX_{/S}$.

Moreover, consider
the composite isomorphism
\begin{equation} \label{z0553}
\eta^\triangledown_{\mfF^\diamondsuit}  :  \BB^{\triangledown, n}  \isom (\mcT_{X^\mr{log}/S^\mr{log}}^{\otimes (n-1)}\otimes  \BB)^\vee  \isom (\mcT_{X^\mr{log}/S^\mr{log}}^{\otimes (n-1)}\otimes \mcF^{n-1})^\vee \isom  (\mcF^0 /\mcF^1)^\vee \isom  (\mcF^\vee)^{n-1} 
\end{equation}
(cf. (\ref{y01}) for the definition of $\BB^{\triangledown, n}$),
where the second isomorphism arise from $\eta_{\mcF^\diamondsuit} :  \BB \isom \mcF^{n-1} $ and the third isomorphism denotes  the dual of  the isomorphism (\ref{d13}) for the case where $j=0$.
One verifies that the collection of data
\begin{equation}
\mfF^{\diamondsuit \triangledown} :=  (\mcF^\vee, \nabla^\vee_\mcF, \{ (\mcF^\vee)^j \}_{j=0}^n, \eta^\triangledown_{\mfF^\diamondsuit})
\end{equation}
forms a $(\mr{GL}_n,  \mbU^\triangledown)$-oper on $\mfX_{/S}$ (cf. (\ref{y02}) for the definition of $\mbU^\triangledown$).
We shall refer to $\mfF^{\diamondsuit \triangledown}$ as the {\it dual $(\mr{GL}_n, \mbU^\triangledown)$-oper of $\mfF^\diamondsuit$}.
If, moreover, $\mfF^\diamondsuit$ is dormant, then  $\mfF^{\diamondsuit \triangledown}$ is immediately   verified   to be dormant.

Finally,  there exists a natural isomorphism
\begin{equation} \label{ww01}
(\mfF^{\diamondsuit \triangledown})^\triangledown \isom \mfF^\diamondsuit
\end{equation}
of $(\mr{GL}_n, \mbU)$-opers.

\vspace{5mm}
\subsection{} \label{z060}

Recall (cf.  ~\cite{Wak5}, \S\,4.4) that  the sheaf of {\it  logarithmic crystalline differential operators} 
(or {\it lcdo}'s for short) 
 on $X^\mr{log}$ over $S^\mr{log}$ is the Zariski sheaf 
\begin{equation} \mcD^{< \infty}_{X^\mr{log}/S^\mr{log}}\end{equation}
 on $X$
   generated, as a sheaf of  rings, by $\mcO_X$ and  
$\mcT_{X^\mr{log}/S^\mr{log}}$
 subject to the following relations:
\begin{itemize}
\item[$\bullet$]
$f_1 \ast f_2 = f_1 \cdot f_2$;
\item[$\bullet$]
$f_1 \ast \xi_1 = f_1 \cdot \xi_1$;
\item[$\bullet$]
$\xi_1 \ast \xi_2 -\xi_2 \ast \xi_1=  [\xi_1, \xi_2] $;
\item[$\bullet$]
$f_1 \ast \xi_1-\xi_1 \ast f_1 = \xi_1 (f_1)$,
\end{itemize}
for local sections $f_1$, $f_2 \in \mcO_X$ and $\xi_1$, $\xi_2 \in \mcT_{X^\mr{log}/S^\mr{log}}$, where $\ast$ denotes the multiplication in $ \mcD^{< \infty}_{X^\mr{log}/S^\mr{log}}$.
In a usual sense, the {\it order $(\geq 0)$} of a given  {\it lcdo}
(i.e., a local section of $\mcD^{< \infty}_{X^\mr{log}/S^\mr{log}}$)
 is well-defined.
Hence, $\mcD^{< \infty}_{X^\mr{log}/S^\mr{log}}$ admits, for each $j \geq 0$, 
 the subsheaf
 \begin{equation}
 \mcD_{X^\mr{log}/S^\mr{log}}^{< j} \ \left(\subseteq  \mcD^{< \infty}_{X^\mr{log}/S^\mr{log}}\right)
 \end{equation}  consisting of {\it lcdo}'s of order $< j$.
$\mcD^{< j}_{X^\mr{log}/S^\mr{log}}$ ($j = 0, 1, 2, \cdots, \infty$) admits two different  structures  of  $\mcO_X$-module --- one as given by left multiplication (where we denote this $\mcO_X$-module by $^l\mcD^{<  j}_{X^\mr{log}/S^\mr{log}}$), and  the other given by right multiplication (where  we denote this $\mcO_X$-module by $^r\mcD^{< j}_{X^\mr{log}/S^\mr{log}}$) ---.
In particular,  we have 
\begin{equation} \label{r4567}
\mcD_{X^\mr{log}/S^\mr{log}}^{<0} = 0 \ \ \text{and} \ \  
{^l\mcD}_{X^\mr{log}/S^\mr{log}}^{<1} = {^r\mcD}_{X^\mr{log}/S^\mr{log}}^{<1} = \mcO_X.
\end{equation}
The set  $\{ \mcD_{X^\mr{log}/S^\mr{log}}^{< j} \}_{j\geq 0}$ forms an increasing filtration on $\mcD^{< \infty}_{X^\mr{log}/S^\mr{log}}$ satisfying that
 \begin{equation} \label{GL7}
 \bigcup_{j \geq 0} \mcD_{X^\mr{log}/S^\mr{log}}^{< j} = \mcD^{< \infty}_{X^\mr{log}/S^\mr{log}},  \ \text{and} \  \ \mcD_{X^\mr{log}/S^\mr{log}}^{< j+1}/\mcD_{X^\mr{log}/S^\mr{log}}^{< j} \isom \mcT_{X^\mr{log}/S^\mr{log}}^{\otimes j}
 \end{equation} 
 for any $j \geq 0$.

Let $\mcF$ be a  vector bundle on $X$, 
 and consider the tensor product ${\mcD}^{<j}_{X^\mr{log}/S^\mr{log}} \otimes \mcF$
  of $\mcF$ with  the $\mcO_X$-module ${^r \mcD}^{<j}_{X^\mr{log}/S^\mr{log}}$.
In the following,  
we shall regard the ${\mcD}^{<j}_{X^\mr{log}/S^\mr{log}} \otimes \mcF$
 as being equipped with a  structure of  $\mcO_X$-module arising from the structure of $\mcO_X$-module ${^l\mcD}^{<j}_{X^\mr{log}/S^\mr{log}}$
   on ${\mcD}^{<j}_{X^\mr{log}/S^\mr{log}}$.

Next, $\nabla_\mcF$ be an $S$-log connection on $\mcF$.
 One may associate $\nabla_\mcF$  with a structure of left ${\mcD}^{<\infty}_{X^\mr{log}/S^\mr{log}}$-module 
 \begin{equation}
 \label{GL70}
 \nabla^\mcD : {\mcD}^{<\infty}_{X^\mr{log}/S^\mr{log}} \otimes \mcF \migi \mcF,
 \end{equation}
(which is $\mcO_X$-linear) on $\mcF$  determined uniquely by the condition that
$\nabla^\mcD (\partial \otimes v) = \langle \nabla (v), \partial \rangle$ for any local sections $v \in \mcF$ and $\partial \in \mcT_{X^\mr{log}/S^\mr{log}}$, where $\langle -,- \rangle$ denotes the
pairing $(\Omega_{X^\mr{log}/S^\mr{log}}\otimes \mcF) \times \mcT_{X^\mr{log}/S^\mr{log}}\migi \mcF$ induced by
the natural paring $\mcT_{X^\mr{log}/S^\mr{log}} \times \Omega_{X^\mr{log}/S^\mr{log}} \migi \mcO_X$.
This assignment $\nabla \mapsto \nabla^\mcD$ determines 
 a  bijective correspondence between the set of $S$-log connections on $\mcF$ and the set of 
  structures of left ${\mcD}^{<\infty}_{X^\mr{log}/S^\mr{log}}$-module 
 $\mcD^{< \infty}_{X^\mr{log}/S^\mr{log} }\otimes \mcF \migi \mcF$ on $\mcF$.


\vspace{5mm}
\subsection{} \label{z061}

Let $\BB$ be an arbitrary line bundle on $X$, and recall  the $p$-determinant data $\mbU_\BB^\mr{can} := (\BB, \nabla^{\mr{can}}_{0, \BB})$ constructed in  \S\,\ref{z018}.
Then, 
 one may construct a canonical $(\mr{GL}_p, \mbU^\mr{can}_\BB)$-oper as follows.
 
First, observe  that there exists  an $\mcO_{X}$-linear morphism 
\begin{align}
\Psi : \mcT^{\otimes p}_{X^\mr{log}/S^\mr{log}}\otimes \BB  \migi \mcD_{X^\mr{log}/S^\mr{log}}^{< \infty}\otimes \BB
\end{align}
determined uniquely by assigning
$\partial^{\otimes p}\otimes b \mapsto (\partial^p - \partial^{(p)}) \otimes b$ for any local sections $\partial \in \mcT_{X^\mr{log}/S^\mr{log}}$, $b \in \BB$.
We shall write
\begin{equation}
\mcD^\Psi_{\BB}
\end{equation}
for 
the quotient of the left $\mcD^{< \infty}_{X^\mr{log}/S^\mr{log}}$-module $ \mcD_{X^\mr{log}/S^\mr{log}}^{< \infty}\otimes \BB$ by the $ \mcD_{X^\mr{log}/S^\mr{log}}^{< \infty}$-submodule generated by the image of $\Psi$.

The structure of left $\mcD^{< \infty}_{X^\mr{log}/S^\mr{log}}$-module on $\mcD^\Psi_{\BB}$
corresponds (cf. \S\,\ref{z060}) to an $S$-log connection
\begin{equation}
\nabla_{\mcD_\BB^\Psi} : \mcD^\Psi_{\BB} \migi \Omega_{X^\mr{log}/S^\mr{log}}\otimes \mcD^\Psi_{\BB}.
\end{equation}

Next, let us write $\mcD^{\Psi, j}_\BB$ for the image of $\mcD_{X^\mr{log}/S^\mr{log}}^{< p-j}\otimes \BB$ via the natural surjection $\mcD_{X^\mr{log}/S^\mr{log}}^{< \infty}\otimes \BB  \migisurj \mcD^\Psi_{\BB}$.
 $\{ \mcD^{\Psi, j}_\BB\}_{j=0}^p$ forms  a $p$-step decreasing filtration 
\begin{equation}
0 =  \mcD^{\Psi, p}_\BB \subseteq \mcD^{\Psi, p-1}_\BB \subseteq \cdots \subseteq \mcD^{\Psi, 0}_\BB = \mcD^{\Psi}_\BB
\end{equation}
on $\mcD^\Psi_{\BB}$ by vector bundles on $X$.

Finally, denote by 
\begin{equation}
\eta_{\mfD^{\Psi \diamondsuit}_\BB} : \BB \isom \mcD^{\Psi, p-1}_\BB 
\end{equation}
the isomorphism
 obtained by restricting the surjection
$\mcD_{X^\mr{log}/S^\mr{log}}^{< \infty}\otimes \BB \migisurj \mcD^\Psi_{\BB}$ to $\mcD^{< 1}_{X^\mr{log}/S^\mr{log}} \otimes \BB$ ($= \BB$) (cf. (\ref{r4567})).

Then,  the following proposition holds.
\vspace{3mm}
\bpr \label{778d3} \leavevmode\\
 \ \ \ 
 The collection of data 
\begin{equation} \label{w01}
  \mfD_\BB^{\Psi \diamondsuit} := (\mcD_\BB^\Psi, \nabla_{\mcD_\BB^\Psi}, \{ \mcD^{\Psi, j}_\BB\}_{j=0}^p, \eta_{\mfD_\BB^{\Psi \diamondsuit}})
  \end{equation}
   forms a dormant $(\mr{GL}_p, \mbU^\mr{can}_\BB)$-oper on $\mfX_{/S}$ of exponent $(\mbF_p, \mbF_p, \cdots, \mbF_p)$.
\epr
\begin{proof}
One may verifies, by the various definitions involved,  that  $\mfD_\BB^{\Psi \diamondsuit}$ is  a dormant $(\mr{GL}_p, \mbU_\BB^\mr{can})$-oper on $\mfX_{/S}$.
The remaining portion (i.e., $\mfD_\BB^{\Psi \diamondsuit}$ is  of exponent $(\mbF_p, \mbF_p, \cdots, \mbF_p)$) follows from Proposition \ref{z047} (since a subset of $\mbF_p$ with cardinality $p$ coincides  with $\mbF_p$ itself).
\end{proof}

\vspace{3mm}
\bpr \label{d3qwe} \leavevmode\\
 \ \ \ 
Let us identify the $p$-determinant data  $\mbU^\mr{can}_{\BB^{\bigstar}}$ with $(\mbU^\mr{can}_\BB)^\triangledown$ via the isomorphism asserted in  Proposition \ref{z03055}.
Then, there exists a canonical isomorphism
\begin{equation} \label{f0f0f}
  \mfD_{\BB^{\bigstar}}^{\Psi \diamondsuit} \isom   (\mfD_{\BB}^{\Psi \diamondsuit})^\triangledown
\end{equation}
of $(\mr{GL}_p, (\mbU^\mr{can}_\BB)^\triangledown)$-opers, where 
\begin{equation}
(\mfD_{\BB^{}}^{\Psi \diamondsuit})^\triangledown := (\mcD_{\BB^{}}^{\Psi \vee}, \nabla_{\mcD_{\BB^{}}^\Psi}^\vee, \{ (\mcD^{\Psi \vee}_{\BB^{}})^j\}_{j=0}^p, \eta^\triangledown_{\mfD_{\BB^{}}^\Psi})
\end{equation}
 denotes the dual $(\mr{GL}_p,  (\mbU^\mr{can}_\BB)^\triangledown)$-oper of $\mfD_{\BB}^{\Psi \diamondsuit}$ (cf. \S\,\ref{z0201}).
 
 \epr
\begin{proof}
Consider the $\mcO_X$-linear  composite
\begin{equation}
\mcD^{< \infty}_{X^\mr{log}/S^\mr{log}} \otimes \BB^{\bigstar}  \isom \mcD^{< \infty}_{X^\mr{log}/S^\mr{log}} \otimes  (\mcD^{\Psi \vee}_{\BB^{}})^{p-1}  \migiincl  \mcD^{< \infty}_{X^\mr{log}/S^\mr{log}} \otimes      \mcD_{\BB^{}}^{\Psi \vee}     \migi \mcD_{\BB^{}}^{\Psi \vee},
\end{equation}
where
\begin{itemize}
\item[$\bullet$]
the first arrow denotes the tensor product of the identity map of $\mcD^{< \infty}_{X^\mr{log}/S^\mr{log}}$ and   the composite isomorphism
\begin{equation}
\BB^{\bigstar} \isom (\mcT_{X^\mr{log}/S^\mr{log}}^{\otimes (p-1)}\otimes \BB^{})^\vee   \isom (\mcD^{\Psi, 0}_{\BB^{}}/\mcD^{\Psi, 1}_{\BB^{}})^\vee \isom  (\mcD^{\Psi \vee}_{\BB^{}})^{p-1}\end{equation}
arising from the various definitions involved;
\item[$\bullet$]
the second arrow denotes the natural inclusion;
\item[$\bullet$]
the third arrow denotes the morphism defining  the structure   of  left $ \mcD^{< \infty}_{X^\mr{log}/S^\mr{log}}$-module  on  $\mcD_{\BB}^{\Psi \vee}$ corresponding to $\nabla_{\mcD_{\BB}^\Psi}^\vee$ (cf. \S\,\ref{z060}).
\end{itemize}
Since $(\mcD_{\BB}^{\Psi \vee}, \nabla_{\mcD_{\BB}^\Psi}^\vee)$ has vanishing $p$-curvature, this composite  turns out to  factor through the quotient $\mcD^{< \infty}_{X^\mr{log}/S^\mr{log}} \otimes \BB^{\bigstar} \migisurj \mcD_{\BB^{\bigstar}}^{\Psi}$.
The resulting morphism
\begin{equation} \label{w0340rt}
\alpha_\BB :  \mcD_{\BB^{\bigstar}}^{\Psi} \migi \mcD_{\BB}^{\Psi \vee},
\end{equation}
 determines an isomorphism of $(\mr{GL}_p, (\mbU^\mr{can}_\BB)^\triangledown)$-opers from $ \mfD_{\BB^{\bigstar}}^{\Psi \diamondsuit}$ to $(\mfD_{\BB}^{\Psi})^\triangledown$, as desired.
\end{proof}

\vspace{6mm}
\section{Duality for dormant $\mr{GL}_n$-opers} \vspace{3mm}

In this section , we discuss duality between   dormant $\mr{GL}_n$-opers and 
 dormant $\mr{GL}_{p-n}$-opers.
Let $S$, $\mfX_{/S}$, $n$, $\mbU$ be as before.
 Suppose further that $n <p$ and $\mbU$ is dormant.
(It follows from Proposition \ref{z03055ww} that such a $\mbU$ necessarily exists.)
\vspace{5mm}
\subsection{} \label{z063}

Let us consider a procedure for constructing a dormant $\mr{GL}_{(p-n)}$-oper  by means of  a dormant $\mr{GL}_n$-oper.

Let  $\vec{\tau} := (\tau_i)_{i=1}^r$ be an element of $2^{\mbF_p}_{\sharp (-) = n}$
 (where $\vec{\tau} := \emptyset$ if $r =0$) and 
$\mfF^\diamondsuit := (\mcF, \nabla_\mcF, \{ \mcF^j \}_{j=0}^n, \eta_{\mfF^\diamondsuit})$  a dormant $(\mr{GL}_n, \mbU)$-oper on $\mfX_{/S}$ of exponent $\vec{\tau}$.
Consider the composite
\begin{equation} \label{d10}
\mcD^{< \infty}_{X^\mr{log}/S^\mr{log}}\otimes \BB \stackrel{\mr{id} \otimes \eta_{\mfF^\diamondsuit}}{\migi} \mcD^{< \infty}_{X^\mr{log}/S^\mr{log}}\otimes \mcF^{n-1} \migiincl \mcD^{< \infty}_{X^\mr{log}/S^\mr{log}}\otimes \mcF \migi \mcF, 
\end{equation}
where the second arrow denotes the natural inclusion and the third arrow
 denotes the morphism defining the structure of left $\mcD^{< \infty}_{X^\mr{log}/S^\mr{log}}$-module on $\mcF$ corresponding to $\nabla_\mcF$.
Since $(\mcF, \nabla_\mcF)$ has vanishing $p$-curvature, this composite turns out to 
 factor through
the  natural surjection $\mcD^{< \infty}_{X^\mr{log}/S^\mr{log}}\otimes \BB \migisurj \mcD^{\Psi}_{\BB}$.
We denote the resulting morphism  by
\begin{equation} \label{z0989}
\widetilde{\eta}_{\mfF^\diamondsuit} : \mcD^{\Psi}_\BB \migi \mcF,
\end{equation}
which is surjective and compatible with the respective $S$-log connections $\nabla_{\mcD^{\Psi}_\BB}$ and $\nabla_\mcF$.
In particular, by restricting $\nabla_{\mcD^{\Psi}_\BB}$, one may construct an $S$-log connection 
\begin{equation}
 \nabla_{\mr{Ker}(\widetilde{\eta}_{\mfF^\diamondsuit} )} : \mr{Ker}(\widetilde{\eta}_{\mfF^\diamondsuit}) \migi \Omega_{X^\mr{log}/S^\mr{log}}\otimes \mr{Ker}(\widetilde{\eta}_{\mfF^\diamondsuit} )
\end{equation}
on   $\mr{Ker}(\widetilde{\eta}_{\mfF^\diamondsuit} )$.
Moreover, for each $j = 0, \cdots, p$, we shall write
\begin{equation}
\mr{Ker}(\widetilde{\eta}_{\mfF^\diamondsuit} )^j :=  \mcD^{\Psi,  j}_\BB \cap \mr{Ker}(\widetilde{\eta}_{\mfF^\diamondsuit} ).
\end{equation}
The inclusions $\mr{Ker}(\widetilde{\eta}_{\mfF^\diamondsuit} )^{j} \migiincl  \mcD^{\Psi,  j}_\BB$ of the   cases where  $j = p-n -1$ and  $p-n$  give rise to a composite
\begin{align} \label{y16}
\eta^{\triangleright}_{\mfF^\diamondsuit} : \mr{Ker}(\widetilde{\eta}_{\mfF^\diamondsuit} )^{(p-n)-1} & \migisurj \mr{Ker}(\widetilde{\eta}_{\mfF^\diamondsuit} )^{(p-n)-1}/\mr{Ker}(\widetilde{\eta}_{\mfF^\diamondsuit} )^{p-n} \\
& \migiincl  \mcD^{\Psi,  (p-n)-1}_\BB /\mcD^{\Psi,  p-n}_\BB \notag  \\
&  \isom \mcD^{\Psi,  p-1}_\BB \otimes \mcT_{X^\mr{log}/S^\mr{log}}^{\otimes n}  \notag \\
& = \BB^{\triangleright, n} \notag 
\end{align}
(cf. (\ref{y03}) for the definition of $ \BB^{\triangleright, n}$),  where the third arrow denotes the  isomorphism (\ref{d13}) of  the case  where the $\mr{GL}_n$-oper is taken to be 
 the underlying  $\mr{GL}_n$-oper  $\mfD_\BB^{\Psi \diamondsuit \heartsuit}$ of $\mfD_\BB^{\Psi \diamondsuit}$.
Thus, we obtain a collection of data
\begin{equation}
\mfF^{\diamondsuit \triangleright}  := (\mr{Ker}(\widetilde{\eta}_{\mfF^\diamondsuit}), \nabla_{\mr{Ker}(\widetilde{\eta}_{\mfF^\diamondsuit})}, \{ \mr{Ker}(\widetilde{\eta}_{\mfF^\diamondsuit})^j \}_{j=0}^{p-n}, \eta^{\triangleright}_{\mfF^\diamondsuit}).
\end{equation}
\vspace{3mm}
\bpr \label{z067}\leavevmode\\
 \ \ \ 
 $\mfF^{\diamondsuit \triangleright}$ is a dormant $(\mr{GL}_{(p-n)}, \mbU^\triangleright)$-oper  on $\mfX_{/S}$ of exponent $\vec{\tau}^{\, \triangleright}$ (cf. (\ref{y04}) for the definition of $\mbU^\triangleright$ and (\ref{w034}) for the definition of $\vec{\tau}^{\, \triangleright}$).
\epr
\begin{proof}
First, let us  prove the claim that
{\it the natural composite
\begin{equation} \label{y12}
\mr{Ker}(\widetilde{\eta}_{\mfF^\diamondsuit}) \migiincl \mcD^\Psi_{\BB} \migisurj  \mcD^\Psi_{\BB}/ \mcD^{\Psi, p-n}_{\BB}
\end{equation}
is an isomorphism}.
To this end, (since both $\mr{Ker}(\widetilde{\eta}_{\mfF^\diamondsuit}) $ and $\mcD^\Psi_{\BB}/ \mcD^{\Psi, p-n}_{\BB}$ are flat over $S$) we may assume, by considering various fibers over $S$, 
that $S = \mr{Spec}(k')$ for  an algebraically closed  field $k'$ over $k$.  
Write
\begin{equation}
\mr{gr}^j := \mr{Ker}(\widetilde{\eta}_{\mfF^\diamondsuit})^j/\mr{Ker}(\widetilde{\eta}_{\mfF^\diamondsuit})^{j+1}
\end{equation}
 ($j= 0, \cdots, p-1$).
The  inclusion $\mr{Ker}(\widetilde{\eta}_{\mfF^\diamondsuit}) \migiincl \mcD^\Psi_\BB$ yields  an inclusion
 \begin{equation} \label{y15}
  \mr{gr}^j \migiincl \mcD^{\Psi, j}_{\BB}/\mcD^{\Psi, j+1}_{\BB}
  \end{equation}
into the line bundle $\mcD^{\Psi, j}_{\BB}/\mcD^{\Psi, j+1}_{\BB}$.
 If $ x_1, \cdots, x_L$ denote the  generic points of irreducible components of $X$, then
the stalk  $\mr{gr}^j_{x_l}$ of $\mr{gr}^j$ at $x_l$ ($l =1, \cdots, L$)  is either trivial or 
free of rank one.
In particular, since the stalk of $\mr{Ker}(\widetilde{\eta}_{\mfF^\diamondsuit})$ at $x_l$ is free of rank $n$, the cardinality of the set $I_l := \big\{ j  \ \big|  \ \mr{gr}^j_{x_l} \neq 0\big\}$ is exactly $n$.
Here, 
recall that the inclusion $\mr{Ker}(\widetilde{\eta}_{\mfF^\diamondsuit}) \migiincl \mcD^\Psi_\BB$ is tautologically compatible with the respective $k'$-log connections $\nabla_{\mr{Ker}(\widetilde{\eta}_{\mfF^\diamondsuit})}$ and $\nabla_{\mcD^\Psi_\BB}$.
Thus, it follows from Proposition \ref{778d3} that $\mr{gr}^{j+1}_{x_l} \neq 0$ implies $\mr{gr}^j_{x_l} \neq 0$.
But this implies that  $I_l= \{ 0,1, \cdots, n-1\}$, and hence, 
 the composite (\ref{y12})
is an isomorphism at every generic point in $X$.
Let us  observe  two sequences of equalities
\begin{align}
& \ \ \ \   \mr{deg}(\mr{Ker}(\widetilde{\eta}_{\mfF^\diamondsuit}) \\
& =  \mr{deg} (\mcD^\Psi_\BB) - \mr{deg} (\mcF) \notag \\
& = p (p-1) (1-g) + p \cdot \mr{deg} (\BB) - n (n-1)(1-g) - n \cdot \mr{deg} (\BB) \notag \\
& = (p-n) (p+n-1 + \mr{deg} (\BB)) \notag
 \end{align}
(where the second equality follows from  (\ref{y13}))
and 
\begin{align}\label{00v}
&   \ \ \ \ \mr{deg}(\mcD_\BB^\Psi/\mcD_\BB^{\Psi, p-n})  \\
& = \mr{deg} (\mcD_\BB^\Psi) - \sum_{j =0}^{n-1}   \mr{deg} (\mcD_\BB^{\Psi, p-j-1} /\mcD_\BB^{\Psi, p -j})\notag \\
& =  \mr{deg} (\mcD_\BB^\Psi) - \sum_{j =0}^{n-1} (\mcT_{X^\mr{log}/S^\mr{log}}^{\otimes j}\otimes \BB) \notag  \\
 & = p (p-1) (1-g) + p \cdot \mr{deg} (\BB) -  \sum_{j =0}^{n-1}   (j (2g-2) + \mr{deg}(\BB) ) \notag  \\
 & =  p (p-1) (1-g) + p \cdot \mr{deg} (\BB) - n(n-1) (g-1) -n \cdot \mr{deg}(\BB) \notag \\
 & =  (p-n) (p+n-1 + \mr{deg} (\BB)). \notag
 \end{align}
By comparing 
$\mr{deg}(\mr{Ker}(\widetilde{\eta}_{\mfF^\diamondsuit}))$ with  $\mr{deg}(\mcD_\BB^\Psi/\mcD_\BB^{\Psi, p-n})$,
we conclude that the  composite (\ref{y12}) is  an isomorphism of $\mcO_X$-modules.  This completes the proof of  the claim.

By virtue of  the claim,    the morphism  (\ref{y15}) for any $j$ and   the composite $\eta^{\triangleright}_{\mcF}$ (cf. (\ref{y16}))  turn out to be   isomorphisms. 
Hence,  (since $\mfD^{\Psi \diamondsuit \heartsuit}_\BB$ is a $\mr{GL}_p$-oper)
the morphism $\mr{gr}^{j+1} \migi \Omega_{X^\mr{log}/S^\mr{log}} \otimes \mr{gr}^{j}$ induced naturally by $\nabla_{\mr{Ker}(\widetilde{\eta}_{\mfF^\diamondsuit})}$ (obtained in the same manner as (\ref{GL2})) is an isomorphism.
Moreover, by the definition of $\nabla_0^{\,\triangleright}$ (cf. \S\,\ref{z029}), the composite isomorphism
\begin{align}
\mr{det}(\mr{Ker}(\widetilde{\eta}_{\mfF^\diamondsuit})) & \isom \mcT_{X^\mr{log}/S^\mr{log}}^{\otimes \frac{(p-n)(p-n-1)}{2}} \otimes (\mr{Ker}(\widetilde{\eta}_{\mfF^\diamondsuit}))^{(p-n)-1})^{\otimes (p-n)} \\
& \hspace{-5mm} \stackrel{\mr{id}\otimes   (\eta^{\triangleright}_{\mfF^\diamondsuit})^{\otimes n}}{\migi}\mcT_{X^\mr{log}/S^\mr{log}}^{\otimes \frac{(p-n)(p-n-1)}{2}} \otimes (\BB^\triangleright)^{\otimes (p-n)} \notag
\end{align}
obtained in the same manner as (\ref{y18}) is compatible with the respective $S$-log connections $\mr{det}(\nabla_{\mr{Ker}(\widetilde{\eta}_{\mfF^\diamondsuit}))})$ and $\nabla_0^\triangleright$.
Finally, observe that since $\mfD^{\Psi\diamondsuit}_\BB$ is dormant (cf. Proposition \ref{778d3}), $\nabla_{\mr{Ker}(\widetilde{\eta}_{\mfF^\diamondsuit})}$ has vanishing $p$-curvature.

Consequently, the collection of data
$\mfF^{\diamondsuit \triangleright} $ forms  a dormant $(\mr{GL}_{(p-n)}, \mbU^\triangleright)$-oper, and we  complete the proof of Proposition \ref{z067}.
\end{proof}

\vspace{5mm}
\subsection{} \label{z06e3}

Consequently, 
by taking account of the discussions in \S\,\ref{z0201} and \S\,\ref{z063},
we obtain  a dormant $(\mr{GL}_{(p-n)}, \mbU^\bigstar)$-oper
\begin{equation} \label{c0c0c}
\mfF^{\diamondsuit \bigstar} := (\mfF^{\diamondsuit \triangleright})^\triangledown
\end{equation}
on $\mfX_{/S}$ of exponent $\vec{\tau}^{\, \bigstar}$.
The assignment $(-)^{\bigstar}$ (i.e.,  $\mfF^{\diamondsuit} \mapsto \mfF^{\diamondsuit \bigstar}$ for each dormant $(\mr{GL}_n,  \mbU)$-oper $\mfF^\diamondsuit$) is compatible with both  base-changing $s'^* (-)$ (cf. (\ref{w0090909}))  via any morphism $s' : S' \migi S$ of $k$-schemes and  tensoring  $(-)_{\otimes \mfL}$ (cf. \S\,\ref{z0301}) with any log integrable line bundle $\mfL := (\mcL, \nabla_\mcL)$.
More precisely, there exist a natural isomorphism 
\begin{equation} \label{w041}
s'^* (\mfF^{\diamondsuit \bigstar}) \isom (s'^*(\mfF^\diamondsuit))^\bigstar 
\end{equation}
of $(\mr{GL}_n, s'^*(\mbU))$-opers 
and
a natural isomorphism
 \begin{equation} \label{w042}
  (\mfF^{\diamondsuit \bigstar})_{\otimes \mfL^\vee} \isom (\mfF^{\diamondsuit}_{\otimes \mfL})^\bigstar
\end{equation}
of $(\mr{GL}_n, \mbU^\bigstar \otimes \mfL^\vee)$-opers, where we identify $\mbU^\bigstar \otimes \mfL^\vee$ with $(\mbU \otimes \mfL)^\bigstar$ via the isomorphism asserted in  Proposition \ref{z0030}.
\vspace{3mm}
\bpr  \label{z07311}\leavevmode\\
 \ \ \ 
Let us identify the $n$-determinant data  $(\mbU^\bigstar)^\bigstar$ with $\mbU$ via the isomorphism asserted in Proposition \ref{z030}.
Then, there exists a natural isomorphism 
\begin{equation}
(\mfF^{\diamondsuit \bigstar})^\bigstar \isom \mfF^\diamondsuit
\end{equation}
of $(\mr{GL}_n, \mbU)$-opers.
\epr
\begin{proof}
For simplicity, we shall write 
\begin{equation}
\mfF^{\diamondsuit \triangleright} =: (\mcE, \nabla_\mcE, \{ \mcE^j \}_{j=0}^n, \eta_{\mfE^\diamondsuit}),
\end{equation}
 which is a dormant $(\mr{GL}_{(p-n)}, \mbU^\triangleright)$-oper on $\mfX_{/S}$.
The natural  inclusion $\mfi \mfn \mfc : \mcE \migiincl \mcD_\BB^\Psi$  induces, via the respective quotients, an isomorphism
$\overline{\mfi \mfn \mfc} : \mcE/ \mcE^1 \isom \mcD_\BB^\Psi/\mcD_\BB^{\Psi, 1}$ (cf. the discussion in the proof of Proposition \ref{z067}, which assert that the morphism (\ref{y15}) is an isomorphism).
One verifies that the composite
\begin{equation}
\BB^\bigstar \isom \Omega_{X^\mr{log}/S^\mr{log}}^{\otimes p}\otimes \BB^\vee \isom (\mcD_\BB^\Psi/\mcD_\BB^{\Psi, 1})^\vee \isom (\mcE/ \mcE^1)^\vee =  (\mcE^\vee)^{(p-n)-1}
\end{equation}
coincides with $\eta_{\mfE^\diamondsuit}^\triangledown$, where
\begin{itemize}
\item[$\bullet$]
 the first arrow denotes the first  isomorphism in  (\ref{y0012});
 \item[$\bullet$]
 the second arrow denotes the dual of the isomorphism (\ref{d13}) of the case where the triple ``$(\mfF^\heartsuit, n, j)$" is taken to be ``$(\mfD_\BB^{\Psi \diamondsuit \heartsuit}, p, 0)$";
 \item[$\bullet$]
the third arrow denotes the dual of $\overline{\mfi \mfn \mfc}$.
\end{itemize}
On the other hand, this composite may be naturally considered  as the restriction to $\mcD^{\Psi, p-1}_{\BB^{\bigstar}}$ of the composite $\mfi \mfn \mfc^\vee \circ \alpha_\BB$ of 
the isomorphism $\alpha_\BB : \mcD^\Psi_{\BB^{\bigstar}} \isom \mcD^{\Psi \vee}_\BB$ asserted in Proposition  \ref{d3qwe}  (cf.  (\ref{w0340rt})) 
 and  the dual $\mfi \mfn \mfc^\vee : \mcD^{\Psi \vee}_\BB \migisurj\mcE^\vee$ of $\mfi \mfn \mfc$.
Thus, it follows from Proposition \ref{d3} that
$\widetilde{\eta}_{\mfE^\diamondsuit}^\triangledown = \mfi \mfn \mfc^\vee \circ \alpha_\BB$.
In particular, by taking the kernels of both $\widetilde{\eta}_{\mfE^\diamondsuit}^\triangledown$ and $\mfi \mfn \mfc^\vee$,
we have a canonical isomorphism
\begin{equation}
 ((\mfF^{\diamondsuit \triangleright})^\triangledown)^\triangleright \isom (\mfF^\diamondsuit)^{\triangledown}
\end{equation}
of $\mr{GL}_n$-opers.
This isomorphism  induces a sequence of isomorphism
\begin{equation}
(\mfF^{\diamondsuit \bigstar})^\bigstar =  (((\mfF^{\diamondsuit \triangleright})^\triangledown)^\triangleright)^\triangledown \isom  ((\mfF^\diamondsuit)^{\triangledown})^\triangledown \isom \mfF^\diamondsuit,
\end{equation}
where the last isomorphism follows from (\ref{ww01}).
The completes the proof of Proposition \ref{z07311}.
\end{proof}

\vspace{5mm}
\subsection{} \label{z072}
We shall write $\mfS \mfe \mft$ for the category of (small) sets.
Also, for a scheme (or, more generally, a stack) $S$
 over $k$, write $\mfS \mfc \mfh_{/S}$ for the category of  relative $S$-schemes.
Let us consider the $\mfS \mfe \mft$-valued functor
\begin{equation}
\mfO \mfp^\ZZZ_{\mr{GL}_n, \mbU, \mfX_{/S}} \ \left(\text{resp.}, \mfO \mfp^\ZZZ_{\mr{GL}_n, \mbU, \tau, \mfX_{/S}}\right) : \mfS \mfc \mfh_{/S} \migi \mfS \mfe \mft
\end{equation}
  on $ \mfS \mfc \mfh_{/S}$ which, to any $S$-schemes $t : T \migi S$, assigns the set of isomorphism classes of dormant $(\mr{GL}_n, t^*(\mbU))$-opers on $\mfX_{/T}$ (resp., dormant $(\mr{GL}_n, t^*(\mbU))$-opers on $\mfX_{T}$ of exponent $\vec{\tau}$) (cf. Remark \ref{y002}).
As explained  later (cf. \S\,\ref{z077}),  both $\mfO \mfp^\ZZZ_{\mr{GL}_n, \mbU, \mfX_{/S}}$ and $\mfO \mfp^\ZZZ_{\mr{GL}_n, \mbU, \tau, \mfX_{/S}}$ may be represented by relative finite  $S$-schemes.
By virtue of Proposition \ref{z047}, 
the relative $S$-scheme $\mfO \mfp^\ZZZ_{\mr{GL}_n, \mbU, \mfX_{/S}}$ decomposes into
the disjoint union 
\begin{equation} \label{w200r}
\mfO \mfp^\ZZZ_{\mr{GL}_n,  \mbU, \mfX_{/S}}  = \coprod_{\vec{\tau} \in (2^{\mbF_p}_{\sharp (-) =n })^{\times r}} \mfO \mfp^\ZZZ_{\mr{GL}_n, \mbU, \vec{\tau}, \mfX_{/S}}.
\end{equation}

\vspace{3mm}
\bt  \label{z073}\leavevmode\\
 \ \ \ 
The assignment  $\mfF^\diamondsuit  \mapsto \mfF^{\diamondsuit \bigstar}$
 constructed in \S\,\ref{z06e3} defines  a canonical isomorphism
\begin{equation} \label{g009o0}
\Theta^\bigstar_{\mr{GL}_n, \mbU, \vec{\tau}, \mfX_{/S}} : \mfO \mfp^\ZZZ_{\mr{GL}_n, \mbU, \vec{\tau}, \mfX_{/S}} \isom  \mfO \mfp^\ZZZ_{\mr{GL}_{(p-n)}, \mbU^\bigstar, \vec{\tau}^\bigstar, \mfX_{/S}} 
\end{equation}
 over $S$ satisfying the equality
\begin{equation}
\Theta^\bigstar_{\mr{GL}_{(p-n)}, \mbU^{\bigstar}, \vec{\tau}^\bigstar, \mfX_{/S}}  \circ \Theta^\bigstar_{\mr{GL}_n, \mbU, \vec{\tau}, \mfX_{/S}}  = \mr{id}
\end{equation}
of automorphisms of $\mfO \mfp^\ZZZ_{\mr{GL}_n,  \mbU, \vec{\tau}, \mbU, \mfX_{/S}}$.
\et
\begin{proof}
The assertion follows from the functoriality of  the assignment $\mfF^\diamondsuit  \mapsto \mfF^{\diamondsuit \bigstar}$ (cf.  (\ref{w041})) and Proposition \ref{z07311}.
\end{proof}

Also, the isomorphism $\Theta^\bigstar_{\mr{GL}_n, \mbU, \vec{\tau}, \mfX_{/S}}$ just obtained  satisfies the following property.

\vspace{3mm}
\bpr\label{z074} \leavevmode\\
 \ \ \ 
Let $\vec{a}$ be an element of $\mbF_p^{\times r}$ (where we take $\vec{a} := \emptyset$ if $r =0$) and
 $\mfL := (\mcL, \nabla_\mcL)$ be a log integrable line bundle on $\mfX_{/S}$ of exponent $\vec{a}$. 
 Then, the assignment $\mfF^\diamondsuit \mapsto \mfF^\diamondsuit_{\otimes \mfL}$ (cf. (\ref{e00908}))
 determines an isomorphism
 \begin{equation}
\Xi_{\mr{GL}_n, \mbU, \vec{\tau}, \mfX_{/S}, \mfL} :  \mfO \mfp_{\mr{GL}_n, \mbU, \vec{\tau}, \mfX_{/S}} \isom \mfO \mfp_{\mr{GL}_n, \mbU\otimes \mfL, \vec{\tau}^{+\vec{a}}, \mfX_{/S}} 
 \end{equation}
over $S$.
Moreover,   the  square diagram
 \begin{equation}
 \begin{CD}
\mfO \mfp_{\mr{GL}_n, \mbU, \vec{\tau}, \mfX_{/S}} @> \Theta^\bigstar_{\mr{GL}_n, \mbU, \vec{\tau}, \mfX_{/S}} >>  \mfO \mfp_{\mr{GL}_{(p-n)}, \mbU^\bigstar, \vec{\tau}^\bigstar, \mfX_{/S}} 
 \\
 @V\Xi_{\mr{GL}_n, \mbU, \vec{\tau}, \mfX_{/S}, \mfL}  VV @VV \Xi_{\mr{GL}_n, \mbU, \vec{\tau}, \mfX_{/S}, \mfL}  V
 \\
\mfO \mfp_{\mr{GL}_n, \mbU\otimes \mfL, \vec{\tau}^{+ \vec{a}}, \mfX_{/S}}  @> \alpha \circ \Theta^\bigstar_{\mr{GL}_n, \mbU, \vec{\tau}^{+ \vec{a}}, \mfX_{/S}} >>   \mfO \mfp_{\mr{GL}_{(p-n)}, \mbU^\bigstar \otimes \mfL^\vee, \vec{\tau}^\bigstar, \mfX_{/S}} 
 \end{CD}
 \end{equation}
is commutative, where 
the lower horizontal arrow denotes the composite of $\Theta^\bigstar_{\mr{GL}_n, \mbU \otimes \mfL, \vec{\tau}^{+ \vec{a}}, \mfX_{/S}} $ and the isomorphism
\begin{equation}
\alpha :  \mfO \mfp^\ZZZ_{\mr{GL}_{(p-n)}, (\mbU  \otimes \mfL)^\bigstar, (\vec{\tau}^{+\vec{a}})^\bigstar, \mfX_{/S}} \isom  \mfO \mfp^\ZZZ_{\mr{GL}_{(p-n)}, \mbU^\bigstar \otimes \mfL^\vee, (\vec{\tau}^\bigstar)^{-\vec{a}}, \mfX_{/S}}
\end{equation}
arising, in an evident fashion, from both the equality $(\vec{\tau}^{+\vec{a}})^\bigstar = (\vec{\tau}^\bigstar)^{-\vec{a}}$ (cf. (\ref{z02ee})) and  the isomorphism  $(\mbU \otimes \mfL)^\bigstar \isom \mbU^\bigstar \otimes \mfL^\vee$ asserted in 
Proposition \ref{z0030}.
\epr
\begin{proof}
The former assertion follows from the discussion in \S\,\ref{z0301}.
The latter assertion follows from the isomorphism  (\ref{w042}).
\end{proof}

\vspace{3mm}
Moreover,  the following  corollary of Theorem \ref{z073} holds.

\vspace{3mm}
\bco  \label{z073ffff}\leavevmode\\
 \ \ \ 
 The structure morphism $\mfO \mfp^\ZZZ_{\mr{GL}_{(p-1)}, \mbU, \mfX_{/S}}  \migi S$ of $\mfO \mfp^\ZZZ_{\mr{GL}_{(p-1)}, \mbU, \mfX_{/S}} $ is an isomorphism.
 That is to say, there exists a unique (up to isomorphism)  dormant $\mr{GL}_{(p-1)}$-oper on $\mfX_{/S}$.
\eco
\begin{proof}

The assertion follows from Theorem \ref{z073} of the case where $n =1$ and  the discussion in Remark \ref{z04eer9g}.
\end{proof}

\vspace{6mm}
\section{Duality for dormant $\mr{GL}_n$-opers on a   proper smooth curve} \label{z0s75}\vspace{3mm}

In this section, we shall consider  the  duality $(-)^\bigstar$ (cf. (\ref{c0c0c})) for dormant $\mr{GL}_n$-opers of  the case where the underlying curve is smooth and there are no marked points.
In this case, one may describe 
the assignment $\mfF^{\diamondsuit} \mapsto \mfF^{\diamondsuit \bigstar}$
in a manner different from that provided in the previous section.

Suppose that $r =0$ and $X/S$ is smooth.
In particular,  the log structures on both $X$ and $S$ defined  in 
\S\,\ref{z04} are trivial (hence, $\Omega_{X^\mr{log}/S^\mr{log}} = \Omega_{X/S}$).

\vspace{5mm}
\subsection{} \label{z07g2}

To begin with, we shall describe the dormant $(\mr{GL}_p, \mbU^\mr{can}_{\mcB})$-oper $\mfD^{\Psi \diamondsuit}_{\mcB}$ (cf. (\ref{w01})) in a manner  different from that provided in \S\,\ref{z061}. 
Let $\mcB$ be a line bundle on $X$.
We shall write
 \begin{align}
 \mcA_{\mcB^\vee} := F^*_{X/S} (F_{X/S*} (\mcB^\vee)).
 \end{align}
Since the relative Frobenius morphism $F_{X/S} : X \migi X_S^{(1)}$ is finite and faithfully flat of degree $p$,
one verifies that $\mcA_{\mcB^\vee}$ is a vector bundle on $X$ of rank $p$.
Denote by
 \begin{align}
 \pi_{\mcB^\vee} :  \mcA_{\BB^\vee} \ \left(= F^*_{X/S}(F_{X/S*}(\BB^\vee))\right) \stackrel{}{\migisurj} \BB^\vee
 \end{align}
 the surjection  determined by 
the adjunction relation ``$F_{X/S}^*(-) \dashv F_{X/S*}(-)$".
Then, $\pi_{\mcB^\vee}$ and the  canonical $S$-connection $\nabla^\mr{can}_{F_{X/S*}(\mcB^\vee)}$ (cf. Remark \ref{z08}) on $\mcA_{\mcB^\vee}$ give rise to a $p$-step decreasing  filtration $\{ \mcA_{\mcB^\vee}^j \}_{j=0}^p$ as follows.
\begin{align}
\mcA_{\BB^\vee}^0  &:= \mcA_{\BB^\vee};\\
\mcA_{\BB^\vee}^1 &:= \mr{Ker}(\mcA_{\BB^\vee} \stackrel{\pi_{\mcB^\vee}}{\migisurj} \BB^\vee); \notag \\
\mcA_{\BB^\vee}^j &:= \mr{Ker}(\mcA_{\BB^\vee}^{j-1} \stackrel{\nabla^{\mr{can}}_{F_{X/S*}(\BB^\vee)}|_{\mcA_{\BB^\vee}^{j-1}}}{\longmigi} \Omega_{X/S} \otimes \mcA_{\BB^\vee} \stackrel{}{\migisurj} \Omega_{X/S}\otimes (\mcA_{\BB^\vee}/\mcA_{\BB^\vee}^{j-1})) \notag
\end{align}
($j= 2, \cdots, p$).
In particular, 
$\pi_{\mcB^\vee}$ induces an isomorphism
\begin{align}
\overline{\pi}_{\mcB^\vee} : \mcA_{\mcB^\vee}/\mcA^1_{\mcB^\vee} \isom \mcB^\vee.
\end{align}
It follows from ~\cite{Wak5},  Proposition 9.3.1  (or,  ~\cite{JP}, Theorem 3.1.6, for the case where $S = \mr{Spec} (\overline{k})$ for an algebraically closed field $\overline{k}$), that
the collection of data
\begin{align}
\mfA^\heartsuit_{\mcB^\vee} := (\mcA_{\mcB^\vee}, \nabla^{\mr{can}}_{F_{X/S*} (\mcB^\vee)}, \{ \mcA^j_{\mcB^\vee} \}_{j=0}^p)
\end{align}
 forms a dormant $\mr{GL}_p$-oper on $\mfX_{/S}$.
We obtain an isomorphism
\begin{align} \label{EEE040}
\Omega_{X/S}^{\otimes (p-1)} \otimes (\mcA_{\mcB^\vee}/\mcA^1_{\mcB^\vee}) \isom \mcA^{p-1}_{\mcB^\vee}
\end{align}
determined by 
 the composite isomorphism (\ref{d13}) of the case where  the triple  $(\mfF^\heartsuit, n, j)$ is taken to be $(\mfA^\heartsuit_{\mcB^\vee}, p, 0)$. 
Let us define  $\eta_{\mfA^\heartsuit_{\mcB^\vee}}$ to be the composite
isomorphism
\begin{align}
\eta_{\mfA^\heartsuit_{\mcB^\vee}} : \left(\BB^{\triangledown, p} =\right) \ \Omega_{X/S}^{\otimes (p-1)} \otimes \mcB^\vee \stackrel{\mr{id} \otimes \overline{\pi}_{\mcB^\vee}^{-1}}{\migi}  \Omega_{X/S}^{\otimes (p-1)} \otimes \mcA_{\mcB^\vee}/\mcA_{\mcB^\vee}^1 \stackrel{(\ref{EEE040})}{\migi} \mcA_{\mcB^\vee}^{p-1}.
\end{align}

Then, the following proposition holds.

\vspace{3mm}
\bpr  \label{z07511}\leavevmode\\
 \ \ \ 
The collection of data
\begin{align}
\mfA^\diamondsuit_{\mcB^\vee} := (\mcA_{\mcB^\vee}, \nabla^{\mr{can}}_{F_{X/S*} (\mcB^\vee)}, \{ \mcA^j_{\mcB^\vee} \}_{j=0}^p, \eta_{\mfA^\heartsuit_{\mcB^\vee}})
\end{align}
forms a dormant  $(\mr{GL}_p, \mbU_\mcB^{\mr{can} \triangledown})$-oper on $\mfX_{/S}$.
Moreover, there exists a canonical isomorphism 
\begin{align}
\gamma_{\mcB} :  \mfD^{\Psi \diamondsuit}_\mcB \isom \mfA^{\diamondsuit \triangledown}_{\mcB^\vee}
\end{align}
of $(\mr{GL}_p, \mbU^{\mr{can}}_\mcB)$-opers.
\epr
\begin{proof}
The former assertion follows immediately from the definition of $\eta_{\mfA^\heartsuit_{\mcB^\vee}}$.
We shall prove the latter assertion.
Consider the composite
\begin{align}
\mcD^{< \infty}_{X^\mr{log}/S^\mr{log}} \otimes \mcB
 & \stackrel{\mr{id} \otimes \overline{\pi}_{\mcB^\vee}^\vee}{\migi} 
 \mcD^{< \infty}_{X^\mr{log}/S^\mr{log}} \otimes (\mcA_{\mcB^\vee}/\mcA^1_{\mcB^\vee})^\vee \\
& \hspace{2.5mm} \migiincl \mcD^{< \infty}_{X^\mr{log}/S^\mr{log}} \otimes \mcA^\vee_{\mcB^\vee} \notag  \\
& \hspace{-3.5mm} \stackrel{\nabla_{F_{X/S*}(\mcB^\vee)}^{\mr{can}\vee \mcD}}{\migi} \mcA^\vee_{\mcB^\vee}, \notag
\end{align}
where 
the second arrow arises from the natural inclusion
$(\mcA_{\mcB^\vee}/\mcA^1_{\mcB^\vee})^\vee \migiincl \mcA^\vee_{\mcB^\vee}$.
Since $\nabla_{F_{X/S*}(\mcB^\vee)}^{\mr{can}\vee}$ has vanishing $p$-curvature,
this composite factors through the quotient $\mcD^{< \infty}_{X^\mr{log}/S^\mr{log}} \otimes \mcB \migisurj \mcD_\mcB^\Psi$.
The resulting morphism
\begin{align} \label{EEE057}
\gamma_\mcB : \mcD^{\Psi}_\mcB \migi \mcA^\vee_{\mcB^\vee}
\end{align}
is surjective and compatible with the respective connections $\nabla_{\mcD^\Psi_\mcB}$
and $\nabla^{\mr{can} \vee}_{F_{X/S*}(\mcB^\vee)}$ (cf. the discussion at the beginning of \S\,\ref{z063}).
Since  the vector bundles $\mcD^{\Psi}_\mcB$ and $\mcA^\vee_{\mcB^\vee}$ have the same rank $p$, 
$\gamma_\mcB$ turns out to be  an isomorphism.
One verifies that
$\gamma_\mcB$ is compatible with the respective filtrations $\{ \mcD^{\Psi, j}_\mcB \}_{j=0}^p$ and $\{ (\mcA^\vee_{\mcB^\vee})^j  \}_{j=0}^p$.
Thus, by the definition of $\eta_{\mfA^\heartsuit_{\mcB^\vee}}$, $\gamma_\mcB$ forms an isomorphism 
$ \mfD^{\Psi \diamondsuit}_\mcB \isom \mfA^{\diamondsuit \triangledown}_{\mcB^\vee}$ of $(\mr{GL}_p, \mbU^{\mr{can}}_\mcB)$-opers.
\end{proof}

\vspace{5mm}
\subsection{} \label{z07fg2}

Now, let us fix 
  an $n$-determinant data $\mbU := (\mcB, \nabla_0)$ (where $\mcB$ is as above)  for $\mfX_{/S}$
and 
   a dormant $(\mr{GL}_n, \mbU)$-oper $\mfF^\diamondsuit := (\mcF, \nabla_\mcF, \{ \mcF^j \}_{j=0}^n, \eta_{\mfF^\diamondsuit})$ on $\mfX_{/S}$. 
Since $\nabla_\mcF^\vee$ has vanishing $p$-curvature,
the morphism
\begin{align} \label{EEE01}
F^*_{X/S}(F_{X/S*} (\mr{Ker} (\nabla_\mcF^\vee))) \migi \mcF^\vee
\end{align}
 corresponding,  via the adjunction relation ``$F^*_{X/S}(-) \dashv F_{X/S*}(-)$",
to the natural inclusion $F_{X/S*} (\mr{Ker}(\nabla_\mcF^\vee)) \migiincl F_{X/S*} (\mcF^\vee)$ is an isomorphism  (cf. ~\cite{Kal}, \S\,5, p.\,190, Theorem 5.1).
Thus, we obtain a composite
\begin{align} \label{EEE03}
\lambda_{\mfF^\diamondsuit} : \mcF^\vee \isom F^*_{X/S}(F_{X/S*} (\mr{Ker} (\nabla_\mcF^\vee))) \migi F^*_{X/S} (F_{X/S*} (\mcB^\vee)) \ \left(= \mcA_{\mcB^\vee}\right),
\end{align}
where the first arrow denotes the inverse of (\ref{EEE01})
and the second arrow denotes the pull-back via $F_{X/S}$ of the composite
\begin{align} \label{jjkk}
F_{X/S*}(\mr{Ker} (\nabla_\mcF^\vee)) \migiincl F_{X/S*} (\mcF^\vee) \migisurj F_{X/S*} (\mcF^{n-1 \vee}) \stackrel{F_{X/S*}(\eta_{\mfF^\diamondsuit}^\vee)}{\migi} F_{X/S*} (\mcB^\vee).
\end{align}

Next, let us write 
\begin{align}
\mcH := \mcA_{\mcB^\vee}/\lambda_{\mfF^\diamondsuit} (\mcF^\vee).
\end{align}
If  $\pi_\mcH : \mcA_{\mcB^\vee} \migisurj \mcH$ denotes the natural surjection, then
  $\mcH$ may be equipped with    a $(p-n)$-step decreasing filtration $\{ \mcH^j \}_{j=0}^{p-n}$   given by 
\begin{align}
\mcH^j :=  \pi_\mcH (\mcA_{\mcB^\vee}^{n+j}) \hspace{3mm} \left(j=0, \cdots, p-n\right).
\end{align}
For each $j \in \{0, \cdots, p-n\}$, $\pi_\mcH$ 
restricts to a morphism
\begin{align}
\pi_{\mcH}^j : \mcA_{\mcB^\vee}^{n+ j} \migi \mcH^{j}.
\end{align}
On the other hand, since $\lambda_{\mfF^\diamondsuit} : \mcF^\vee \migi \mcA_{\mcB^\vee}$ is compatible with the respective $S$-connections $\nabla_\mcF^\vee$ and 
$\nabla^\mr{can}_{F_{X/S*}(\mcB^\vee)}$,
$\nabla^\mr{can}_{F_{X/S*}(\mcB^\vee)}$ induces, via $\pi_\mcH$,
 an $S$-connection
 \begin{align}
 \nabla_\mcH : \mcH \migi \Omega_{X/S} \otimes \mcH
 \end{align}
 on $\mcH$.
Then, the following lemma holds.

\vspace{3mm}
\ble\label{zF07L4} \leavevmode\\
 \ \ \ 
For each $j \in \{0, \cdots, p-n\}$, the morphism $\pi^j_{\mcH}$
  is an isomorphism.
Moreover, the collection of data
 \begin{align}
 \mfH^\heartsuit := (\mcH, \nabla_\mcH, \{ \mcH^j \}_{j=0}^{p-n})
 \end{align}
 forms a dormant  $\mr{GL}_{(p-n)}$-oper on $\mfX_{/S}$.
\ele
\begin{proof}
Let us observe the following two facts (a), (b):
\begin{itemize}
\item[(a)]
The square diagram
\begin{align} \label{EEE060}
\begin{CD}
\mcF^\vee @> \lambda_{\mfF^\diamondsuit} >> \mcA_{\mcB^\vee}
\\
@VVV @VV \pi_{\mcB^\vee} V
\\
\mcF^{n-1 \vee} @>> \eta_{\mfF^\diamondsuit}^\vee  > \mcB^\vee
\end{CD}
\end{align}
is commutative, where  the left-hand vertical arrow denotes the dual of the inclusion $\mcF^{n-1} \migiincl \mcF$.
Moreover, the kernels of the  left-hand and right-hand vertical  arrows coincide with $(\mcF^{\vee})^1$ and $\mcA^1_{\mcB^\vee}$ respectively.
In particular, (since $\eta_{\mfF^\diamondsuit}^\vee$ is an isomorphism) the diagram (\ref{EEE060}) induces  an isomorphism
\begin{align} \label{EEE061}
\mcF^\vee / (\mcF^\vee)^1 \isom \mcA_{\mcB^\vee} / \mcA_{\mcB^\vee}^1.
\end{align}
\item[(b)]
For any $j =2, \cdots, n$, (by the definition of a $\mr{GL}_n$-oper) the following equalities holds:
\begin{align}
(\mcF^\vee)^j & = \mr{Ker} ((\mcF^\vee)^{j-1} \stackrel{\nabla_\mcF^\vee |_{(\mcF^\vee)^{j-1}}}{\migi} \Omega_{X/S}\otimes \mcF^\vee \migisurj \Omega_{X/S} \otimes (\mcF^\vee/(\mcF^\vee)^{j-1})), \\
\mcA_{\mcB^\vee}^j & = \mr{Ker} (\mcA_{\mcB^\vee}^{j-1} \stackrel{\nabla_{F_{X/S*}(\mcB^\vee)}^{\mr{can}}|_{\mcA_{\mcB^\vee}^{j-1}}}{\migi} \Omega_{X/S}\otimes \mcA_{\mcB^\vee} \migisurj \Omega_{X/S} \otimes (\mcA_{\mcB^\vee}/\mcA_{\mcB^\vee}^{j-1})). \notag
\end{align}
\end{itemize}
These facts imply, by  induction on $j$,  that the morphism $\lambda_{\mfF^\diamondsuit}$ restricts to morphisms
\begin{align}
(\mcF^\vee)^j \migi \mcA_{\mcB^\vee}^j \hspace{3mm} \left(j =0, \cdots, n\right).
\end{align}
The resulting morphisms 
\begin{align}
\overline{\lambda}_{\mfF^\diamondsuit}^j : (\mcF^\vee)^j /(\mcF^\vee)^{j+1} \migiincl \mcA^j_{\mcB^\vee} / \mcA^{j+1}_{\mcB^\vee} \hspace{3mm} \left(j=0,\cdots, n-1\right)
\end{align}
between the respective subquotients of $\mcF^\vee$ and $\mcA_{\mcB^\vee}$ make the square diagrams
\begin{align} \label{EEE063}
\begin{CD}
(\mcF^\vee)^j /(\mcF^\vee)^{j+1} @> \overline{\lambda}_{\mfF^\diamondsuit}^j >> \mcA^j_{\mcB^\vee} / \mcA^{j+1}_{\mcB^\vee}  
\\
@V \mfk \mfs_{\mfF^\heartsuit}^j V \wr V @V \wr V  \mfk \mfs_{\mfA_{\mcB^\vee}^\heartsuit}^j  V
\\
\Omega_{X/S} \otimes ((\mcF^\vee)^{j-1} /(\mcF^\vee)^{j}) @> \mr{id} \otimes  \overline{\lambda}_{\mfF^\diamondsuit}^{j-1} >> \Omega_{X/S} \otimes (\mcA^{j-1}_{\mcB^\vee} / \mcA^{j}_{\mcB^\vee}) 
\end{CD}
\end{align}
$(j =0, \cdots, n)$ commute.
By induction on $j$, the commutativity of  (\ref{EEE063}) and the isomorphism (\ref{EEE061}) imply that
all the $\overline{\lambda}_{\mfF^\diamondsuit}^j$'s are isomorphisms.
Hence, the composite $\mcF^\vee \stackrel{\lambda_{\mfF^\diamondsuit}}{\migi} \mcA_{\mcB^\vee} \migisurj \mcA_{\mcB^\vee} / \mcA^n_{\mcB^\vee}$ is an isomorphism, equivalently, $\pi^0_{\mcH}$ is an isomorphism.
It follows that all the  $\pi^j_\mcH$'s are  isomorphisms.
This completes the proof of the former assertion.

Next, let us consider the latter assertion.
The morphisms $\pi^j_\mcH$ induce  commutative square diagrams
\begin{align} \label{EEE072}
\begin{CD}
\mcA^{n+j}_{\mcB^\vee} / \mcA^{n+j+1}_{\mcB^\vee}  @>>> \mcH^j /\mcH^{j+1}
\\
@V \mfk \mfs_{\mfA_{\mcB^\vee}^\heartsuit}^j V V @V V \mfk \mfs^j_{\mfH^\heartsuit}V
\\
 \Omega_{X/S} \otimes (\mcA^{n+j-1}_{\mcB^\vee} / \mcA^{n+ j}_{\mcB^\vee})  @>>> \Omega_{X/S} \otimes (\mcH^{j-1}/\mcH^j)
\end{CD}
\end{align}
$(j = 0, \cdots, n-p)$, 
where the upper and lower horizontal arrows are isomorphisms arising from $\pi^j_{\mcH}$ and $\pi^{j-1}_\mcH$ respectively.
Thus, the latter assertion follows from the commutativity of (\ref{EEE072}) and 
 the fact  (cf. Proposition \ref{z07511}) that $(\mcA_{\mcB^\vee}, \nabla^\mr{can}_{F_{X/S*}(\mcB^\vee)}, \{ \mcA^j_{\mcB^\vee} \}_{j=0}^p)$ is a $\mr{GL}_p$-oper (which implies that the $\mfk \mfs_{\mfA_{\mcB^\vee}^\heartsuit}^j$'s are isomorphisms).
\end{proof}

 Finally, let us define  $\eta_{\mfH^\diamondsuit}$ to be the composite isomorphism
 \begin{align}
 \eta_{\mfH^\diamondsuit} :
  \mcB^\bigstar \stackrel{(\ref{y0012})}{\migi}
     \Omega_{X/S}^{\otimes (p-1)}\otimes \mcB^\vee
 \stackrel{\mr{id} \otimes \overline{\pi}_{\mcB^\vee}^{-1}}{\migi} \Omega_{X/S}^{\otimes (p-1)}\otimes (\mcA_{\mcB^\vee}/ \mcA_{\mcB^\vee}^1) 
 \stackrel{(\ref{EEE040})}{\migi} \mcA_{\mcB^\vee}^{p-1} \stackrel{\pi_\mcH^{n-p-1}}{\migi} \mcH^{n-p-1}.
 \end{align}
Then, by the following proposition,
the assignment $\mfF^\diamondsuit \mapsto \mfF^{\diamondsuit \bigstar}$ may be identified with
the assignment $\mfF^\diamondsuit \mapsto \mfH^\diamondsuit$ (cf. (\ref{EEE046})).

\vspace{3mm}
\bpr\label{zF074} \leavevmode\\
 \ \ \ 
The collection of data
\begin{align} \label{EEE046}
\mfH^\diamondsuit := (\mcH, \nabla_\mcH, \{ \mcH^j \}_{j=0}^{p-n}, \eta_{\mfH^\diamondsuit})
\end{align}
forms a dormant $(\mr{GL}_{p-n}, \mbU^\bigstar)$-oper on $\mfX_{/S}$ and  is isomorphic to $\mfF^{\diamondsuit \bigstar}$.
\epr
\begin{proof}
Denote by   $\iota_{\mfD^{\Psi \diamondsuit}_\mcB} : \mcB \migiincl \mcD^{\Psi}_{\mcB}$ (resp., $\iota_{\mfF^\diamondsuit} : \mcB \migiincl  \mcF$)  the injection defined to be
the composite of $\eta_{\mcD_\mcB^{\Psi \diamondsuit}} : \mcB \isom \mcD^{\Psi, p-1}_\mcB$
(resp., $\eta_{\mfF^\diamondsuit} : \mcB \isom \mcF^{n-1}$) and the inclusion 
$\mcD^{\Psi, p-1}_\mcB \migiincl \mcD^{\Psi}_\mcB$ (resp., $\mcF^{n-1} \migiincl \mcF$).
The following equalities hold:
\begin{align}
\widetilde{\eta}_{\mfF^\diamondsuit} \circ \iota_{\mfD^{\Psi \diamondsuit}_\mcB} =\iota_{\mfF^\diamondsuit} =  \lambda^\vee_{\mfF^\diamondsuit}  \circ  \pi^\vee_{\mcB^\vee} = \lambda^\vee_{\mfF^\diamondsuit} \circ\gamma_\mcB \circ  \iota_{\mfD^{\Psi \diamondsuit}_\mcB}. 
\end{align}
Hence, it follows from Proposition  \ref{d3} that
(since both $\widetilde{\eta}_{\mfF^\diamondsuit}$ and $ \lambda^\vee_{\mfF^\diamondsuit} \circ\gamma$ are compatible with the respective connections $\nabla_{\mcD^{\Psi}_\mcB}$ and $\nabla_\mcF$)
 the diagram
\begin{align} \label{EEE011}
\xymatrix{
  \mcD^\Psi_\mcB \ar[rd]_{\widetilde{\eta}_{\mfF^\diamondsuit}} \ar[rr]^{\gamma_\mcB}_{\sim}& &\mcA^\vee_{\mcB^\vee}  \ar[ld]^{\lambda^\vee_{\mfF^\diamondsuit}} \\
& \mcF &
}
\end{align}
is commutative.
This commutative diagram  induces
an isomorphism 
\begin{align} \label{EEE021}
 (\mr{Ker}(\widetilde{\eta}_{\mfF^\diamondsuit}), \nabla_{\widetilde{\eta}_{\mfF^\diamondsuit}}) \isom (\mcH^\vee, \nabla_\mcH^\vee)
\end{align}
 of (log) integrable vector bundles.
It follows from  the various definitions involved that the dual isomorphism 
\begin{align} \label{EEE025}
\delta : (\mr{Ker}(\widetilde{\eta}_{\mfF^\diamondsuit})^\vee, \nabla^\vee_{\widetilde{\eta}_{\mfF^\diamondsuit}}) \isom (\mcH, \nabla_\mcH)
 \end{align}
  of  (\ref{EEE021})
  is compatible with the respective filtrations $\{ (\mr{Ker}(\widetilde{\eta}_{\mfF^\diamondsuit})^\vee)^j \}_{j=0}^{p-n}$ and $\{ \mcH^j \}_{j=0}^{n-p}$, and moreover, satisfies the equality $\delta |_{(\mr{Ker}(\widetilde{\eta}_{\mfF^\diamondsuit})^\vee)^{p-n-1}} \circ \eta^{\triangleright \triangledown}_{\mfF^\diamondsuit} = \eta_{\mfH^\diamondsuit}$.
  This completes the proof of Proposition \ref{zF074}.
  \end{proof}

\vspace{3mm}

\begin{rema} \label{z0hhy49} \leavevmode\\
 \ \ \ 
Once  the  proposition just above is  proved, one may describe the assignment $\mfF^\diamondsuit \mapsto \mfH^\diamondsuit$ ($\cong \mfF^{\diamondsuit \bigstar}$) more briefly, as follows. 
For simplicity, suppose further  that $S = \mr{Spec} (\overline{k})$ for an algebraically closed field $\overline{k}$ over $k$.
Let $\mbU$ and $\mfF^\diamondsuit$ be as above.
Denote by $\alpha  : F_{X/\overline{k}*}(\mr{Ker} (\nabla_\mcF^\vee)) \migi F_{X/\overline{k}*} (\mcB^\vee)$ the composite
(\ref{jjkk}).
Then, 
one verifies that $\mfH^\diamondsuit$ is isomorphic to
\begin{align}
(F^*_{X/\overline{k}}(\mr{Coker} (\alpha)), \nabla^\mr{can}_{\mr{Coker} (\alpha)}, \{ F^*_{X/\overline{k}}(\mr{Coker} (\alpha))^{[j]} \}_{j=0}^{p-n}),
\end{align}
where $\{ F^*_{X/\overline{k}}(\mr{Coker} (\alpha))^{[j]} \}_{j=0}^{p-n}$
denotes  the 
 Harder-Narasimhan filtration on the vector bundle  $F^*_{X/\overline{k}}(\mr{Coker} (\alpha))$.
\end{rema}

\vspace{6mm}
\section{Duality for dormant $\mfs \mfl_n$-opers} \label{z075}\vspace{3mm}

In ~\cite{Wak5}, we discussed the definition of a(n) (dormant)  $\mfs \mfl_n$-oper (or, more generally, a $\mfg$-oper for a semisimple Lie algebra $\mfg$) on a pointed stable curve and various properties concerning their moduli.
In this last section, we consider duality for  dormant $\mfs \mfl_n$-opers  induced by  Theorem \ref{z073} and 
some applications (by means of results of $p$-adic Teichm\"{u}ller theory) to understanding the moduli stack of dormant $\mfs \mfl_{n}$-opers of the case where $n = p-2$.
\vspace{5mm}
\subsection{} \label{z076}
Suppose that $p-1  > n>1$.
We shall identify the Lie algebra $\mfp \mfg \mfl_n$ of $\mr{PGL}_n$  with 
the spectrum of the symmetric algebra $\mbS_k (\mfp \mfg \mfl_n^\vee)$ on $\mfp \mfg \mfl_n$ over $k$.
Since $n <p$,  one may  identify $\mfs \mfl_n$ with $\mfp \mfg \mfl_n$ via the natural composite  $\mfs \mfl_n \migiincl \mfg \mfl_n \migisurj   \mfp \mfg \mfl_n$ and moreover,  identify $\mr{PGL}_n$ with the adjoint group of $\mfs \mfl_n$.
Write $\mfc_n$ for  the GIT quotient $\mfp \mfg \mfl_n \ooalign{$/$ \cr $\,/$}\hspace{-0.5mm}\mr{PGL}_n$ of $\mfp \mfg \mfl_n$ by the adjoint action of $\mr{PGL}_n$, i.e., the spectrum 
$\mr{Spec}(\mbS_k (\mfp \mfg \mfl_n^\vee)^{\mr{PGL}_n})$
of the ring of polynomial invariants on $\mfp \mfg \mfl_n$.
Denote by 
\begin{equation}
\chi_n : \mfp \mfg \mfl_n \migi \mfc_{n}
\end{equation}
the so-called {\it Chevalley map}, i.e.,  the morphism over $k$ corresponding to the inclusion $\mbS_k (\mfp \mfg \mfl_n^\vee)^{\mr{PGL}_n} \migiincl \mbS_k (\mfp \mfg \mfl_n^\vee)$ of $k$-algebras.
Also, denote by $\mft_{n}$ the Lie subalgebra of $\mfp \mfg \mfl_n$ consisting of  the image (via $\mfg \mfl_n \migisurj \mfp \mfg \mfl_n$) of  diagonal matrices, and by $\kappa : \mft_n \migiincl \mfp \mfg \mfl_n$ the natural inclusion.

Note that  the various  Lie algebras and morphisms between them mentioned above may be defined over $\mbF_p$. 
Hence, 
it makes sense to speak of the sets of the  $\mbF_p$-rational points $\mfp \mfg \mfl_n (\mbF_p)$, $\mfc_n (\mbF_p)$,  $\mft_n (\mbF_p)$ of $\mfp \mfg \mfl_n$, $\mfc_n$, $\mft_n$ respectively,  and of the maps between the respective   sets of $\mbF_p$-rational points $\chi_n (\mbF_p) : \mfp \mfg \mfl_n (\mbF_p) \migi \mfc_n (\mbF_p)$, $\kappa (\mbF_p) : \mft_n (\mbF_p) \migi \mfp \mfg \mfl_n (\mbF_p)$.
The set $\mft_n (\mbF_p)$
may be identified with 
the quotient   $\mbF_p^{\times n}/\varDelta (\mbF_p)$, where $\varDelta$ denotes the diagonal  embedding $\mbF_p \migiincl \mbF_p^{\times n}$.
The symmetric group $\mfS_n$ of $n$ letters acts, by permutation, on  $\mbF_p^{\times n}/\varDelta (\mbF_p)$ in such a way that the surjection $\mbF_p^{\times n} \migisurj \mbF_p^{\times n}/\varDelta (\mbF_p)$ is compatible with the respective $\mfS_n$-actions (cf. (\ref{w101})).
Hence, we have the double quotient $\mfS_n \backslash \mbF_p^{\times n}/\varDelta (\mbF_p)$.
The natural surjection $\mbF_p^{\times n} \migisurj \mfS_n \backslash \mbF_p^{\times n}/\varDelta (\mbF_p)$ factors through the surjection $\mbF_p^{\times n} \migisurj \mbN_{\sharp ( -)=n}^{\mbF_p}$   (cf. (\ref{w101})), and we have the resulting surjection 
\begin{equation} \label{w104}
\pi^{\mbF_p}_{n} : \mbN_{\sharp ( -)=n}^{\mbF_p} \migisurj  \mfS_n \backslash \mbF_p^{\times n}/\varDelta (\mbF_p).
\end{equation}


Next, observe that the composite
\begin{equation} \label{w103}
\mbF_p^{\times n}/\varDelta (\mbF_p) \ \left(=\mft_{n} (\mbF_p)\right) \stackrel{\kappa (\mbF_p)}{\migiincl}  \mfp \mfg \mfl_n (\mbF_p) \stackrel{\chi_n (\mbF_p)}{\migi} \mfc_{n} (\mbF_p)
\end{equation}
factors through the surjection  $\mbF_p^{\times n}/\varDelta (\mbF_p) \migisurj \mfS_n \backslash \mbF_p^{\times n}/\varDelta (\mbF_p)$.
The resulting map of sets
\begin{equation} \label{w102}
\chi^{\mbF_p}_n : \mfS_n \backslash \mbF_p^{\times n}/\varDelta (\mbF_p) \migi \mfc_{n} (\mbF_p)
\end{equation}
is injective.

For each element $a  \in \mbF_p$ and $\rho_0  \in \mr{Im}(\chi_n^{\mbF_p})$ ($\subseteq \mfc_{n} (\mbF_p)$), there exists 
a {\it unique}  multiset 
\begin{equation} \label{v01}
\rho_0^{\divideontimes a } := [\rho_{0, 1}^{\divideontimes a }, \cdots, \rho_{0, n}^{\divideontimes a }]
\end{equation}
(where $\rho_{0, j}^{\divideontimes a } \in \mbF_p$ for each $j \in \{ 1, \cdots, n \}$)
over $\mbF_p$ with cardinality $n$ (i.e.,
an element of $ \mbN_{\sharp (-)=n}^{\mbF_p}$)  satisfying  that $\chi_n^{\mbF_p}\circ \pi_n^{\mbF_p}(\rho_0^{\divideontimes a }) = \rho$ and $\sum_{j=1}^n \rho_{0, j}^{\divideontimes a } = a$. 

Let us define a subset $\mfc_{n}(\mbF_p)^{\circledast}$ of $\mfc_{n} (\mbF_p)$ to be
\begin{equation} \label{v0032}
\mfc_{n}(\mbF_p)^{\circledast} : =\chi_n^{\mbF_p}\circ \pi_n^{\mbF_p} (2^{\mbF_p}_{\sharp (-) = n}).
\end{equation}
If $\rho_0 \in \mfc_{n}(\mbF_p)^{\circledast}$ and $a \in \mbF_p$, 
then (since $\rho_0^{\divideontimes a }$ is a subset of $\mbF_p$) the element
\begin{equation} \label{v0033}
\rho_0^\bigstar :=  \chi_{(p-n)}^{\mbF_p}\circ \pi_{(p-n)}^{\mbF_p}((\rho_0^{\divideontimes a })^\bigstar)
\end{equation}
is well-defined and lies in  $\mfc_{(p-n)}(\mbF_p)^\circledast$.
Note that the element $\rho^\bigstar$ does not depend on the choice of $a$.
Thus, the assignment $\rho \mapsto \rho^\bigstar$ defines a well-defined bijection of sets
\begin{equation}
\theta_n^\bigstar : \mfc_{n}(\mbF_p)^{\circledast} \migi \mfc_{(p-n)}(\mbF_p)^\circledast,
\end{equation}
which is verified to satisfy the equality  $\theta_{(p-n)}^\bigstar \circ \theta_n^\bigstar = \mr{id}$.

For each $\vec{\rho} := (\rho_i)_{i=1}^r\in \mfc_n(\mbF_p)^{\times r}$
and $\vec{a} := (a_i)_{i =1}^r \in \mbF_p^{\times r}$,  we shall write
\begin{equation} \label{v0034}
\vec{\rho}^{\, \divideontimes \vec{a}} := (\rho_{i}^{\divideontimes a_i})_{i=1}^r \in (\mbN_{\sharp ( -)=n}^{\mbF_p})^{\times r}.
\end{equation} 
If, moreover, $\vec{\rho}$ lies in $(\mfc_n(\mbF_p)^\circledast)^{\times r}$, then one obtains
\begin{equation} \label{h8989}
\vec{\rho}^{\, \bigstar}  := (\rho_i^\bigstar)_{i=1}^r \in  (\mfc_{(p-n)}(\mbF_p)^\circledast)^{\times r}.
\end{equation}
In particular, $(\vec{\rho}^{\, \bigstar})^\bigstar = \vec{\rho}$.
\vspace{5mm}
\subsection{} \label{z077}

Let $\vec{\rho} := (\rho_i)_{i=1}^r \in \mfc_n (\mbF_p)^{\times r}$ (where $\vec{\rho} := \emptyset$ if $r = 0$).
We shall write
\begin{equation}
\mfO \mfp^\ZZZ_{\mfs \mfl_n, \mfX_{/S}} \  \left(\text{resp.,} \ \mfO \mfp^\ZZZ_{\mfs \mfl_n, \vec{\rho}, \mfX_{/S}}\right) : \mfS \mfc \mfh_{/S} \migi \mfS \mfe \mft
\end{equation}
for   the moduli functor, as  introduced in ~\cite{Wak5} \S\,3.6, classifying (the isomorphism classes of) dormant $\mfs \mfl_n$-opers on $\mfX_{/S}$  (resp., dormant $\mfs \mfl_n$-opers on $\mfX_{/S}$ of radii $\vec{\rho}$).
(For convenience, we shall say that  any dormant $\mfs \mfl_n$-oper is {\it of radii $\emptyset$}.)
 Both 
$\mfO \mfp^\ZZZ_{\mfs \mfl_n, \mfX_{/S}}$ and $\mfO \mfp^\ZZZ_{\mfs \mfl_n, \vec{\rho}, \mfX_{/S}}$
 may be represented by  (possibly empty)
relative finite $S$-schemes (cf. ~\cite{Wak5}, Theorem C).
Moreover,   
if $r >0$, 
then 
$\mfO \mfp^\ZZZ_{\mfs \mfl_n, \mfX_{/S}}$ decomposes into the disjoint union
\begin{equation} \label{w4405}
\mfO \mfp^\ZZZ_{\mfs \mfl_n, \mfX_{/S}}  = \coprod_{\vec{\rho} \in \mfc_n (\mbF_p)^{\times r}} \mfO \mfp^\ZZZ_{\mfs \mfl_n, \vec{\rho}, \mfX_{/S}}. 
\end{equation}
(cf.  ~\cite{Wak5},  Theorem C  (i)).


Now, let 
$\vec{a} := (a_i)_{i=1}^r \in \mbF_p^{\times r}$, $\vec{\rho} :=(\rho_i)_{i=1}^r \in \mfc_n (\mbF_p)^{\times r}$ (where we take $\vec{a} := \emptyset$ and $\vec{\rho} := \emptyset$ if $r =0$), and let
$\mbU := (\BB, \nabla_0)$ be  a dormant  $n$-determinant data for $\mfX_{/S}$ of exponent $\vec{a}$.
Also, let  $\mfF^\diamondsuit$
 be a dormant $(\mr{GL}_n, \mbU)$-oper on $\mfX_{/S}$ (resp., a dormant $(\mr{GL}_n, \mbU)$-oper on $\mfX_{/S}$ of exponent $\vec{\rho}^{\, \divideontimes \vec{a}}$, where $\vec{\rho} \in \mfc_n(\mbF_p)^{\times r}$).
 It induces, via the  change of structure group $\mr{GL}_n \migisurj \mr{PGL}_n$,  a dormant  $\mfs \mfl_n$-oper on $\mfX_{/S}$ (resp.,  a dormant  $\mfs \mfl_n$-oper on $\mfX_{/S}$ of radii $\vec{\rho}$), which we denote by  $\mcF^{\diamondsuit \spadesuit}$.
The assignment $\mcF^\diamondsuit \mapsto \mcF^{\diamondsuit \spadesuit}$ defines a morphism
\begin{align} \label{w5001}
& \hspace{17mm} \Lambda_{\mr{GL}_n, \mbU, \mfX_{/S}} : \mfO \mfp^\ZZZ_{\mr{GL}_n,  \mbU, \mfX_{/S}} \migi \mfO \mfp^\ZZZ_{\mfs \mfl_n, \mfX_{/S}}  \\
& \left(\text{resp.,} \ \Lambda_{\mr{GL}_n, \mbU, \vec{\rho},  \mfX_{/S}} : \mfO \mfp^\ZZZ_{\mr{GL}_n, \mbU,  \vec{\rho}^{\, \divideontimes \vec{a}}, \mfX_{/S}} \migi \mfO \mfp^\ZZZ_{\mfs \mfl_n, \vec{\rho}, \mfX_{/S}} \right) \notag
\end{align}
over $S$.
By ~\cite{Wak5}, Corollary 4.14.3 (i),  $\Lambda_{\mr{GL}_n, \mbU, \mfX_{/S}} $ is an isomorphism.
(In particular, $\mfO \mfp^\ZZZ_{\mr{GL}_n,  \mbU, \mfX_{/S}}$, as well as $\mfO \mfp^\ZZZ_{\mr{GL}_n, \mbU,  \vec{\rho}^{\, \divideontimes \vec{a}}, \mfX_{/S}}$, may be represented by a relative finite $S$-scheme.)
On the other hand, it follows from  Proposition \ref{z047} that if $r >0$, then $\mfO \mfp^\ZZZ_{\mr{GL}_n, \mbU, \mfX_{/S}}$ decomposes into the disjoint union
\begin{equation}
\mfO \mfp^\ZZZ_{\mr{GL}_n, \mbU, \mfX_{/S}}  = \coprod_{\vec{\rho} \in (\mfc_n (\mbF_p)^\circledast)^{\times r} } \mfO \mfp^\ZZZ_{\mr{GL}_n, \mbU, \vec{\rho}^{\, \divideontimes \vec{a}}, \mfX_{/S}}.
\end{equation}
Hence, the decomposition (\ref{w4405}) may be described as 
\begin{equation} \label{w4660}
\mfO \mfp^\ZZZ_{\mfs \mfl_n, \mfX_{/S}}  = \coprod_{\vec{\rho} \in (\mfc_n(\mbF_p)^\circledast)^{\times r} } \mfO \mfp^\ZZZ_{\mfs \mfl_n, \vec{\rho}, \mfX_{/S}},
\end{equation}
and, by restricting $\Lambda_{\mr{GL}_n, \mbU, \mfX_{/S}}$, we obtain an isomorphism 
\begin{equation} \label{w4661}
\Lambda_{\mr{GL}_n, \mbU, \vec{\rho}, \mfX_{/S}} :  \mfO \mfp^\ZZZ_{\mr{GL}_n, \mbU, \vec{\rho}^{\, \divideontimes \vec{a}}, \mfX_{/S}} \isom  \mfO \mfp^\ZZZ_{\mfs \mfl_n, \vec{\rho}, \mfX_{/S}}
\end{equation}
over $S$.
Consider  the composite isomorphisms 
\begin{align} \label{w4662}
\Theta_{\mfs \mfl_n, \mfX_{/S}}^\bigstar &: = \Lambda_{\mr{GL}_{(p-n)}, \mbU, \mfX_{/S}} \circ \Theta_{\mr{GL}_n, \mbU, \mfX_{/S}}^\bigstar \circ \Lambda_{\mr{GL}_n, \mbU, \mfX_{/S}}^{-1}, \\
\Theta_{\mfs \mfl_n, \vec{\rho}, \mfX_{/S}}^\bigstar &: = \Lambda_{\mr{GL}_{(p-n)}, \mbU, (\vec{\rho}^{\, \divideontimes \vec{a}})^\bigstar, \mfX_{/S}} \circ \Theta_{\mr{GL}_n, \mbU, \vec{\rho}^{\, \divideontimes \vec{a}}, \mfX_{/S}}^\bigstar \circ \Lambda_{\mr{GL}_n, \mbU, \vec{\rho}, \mfX_{/S}}^{-1}. \notag
\end{align}
 Proposition \ref{z074} implies that  these morphisms do not depend on the choice of $\mbU$.
By 
  Theorem \ref{z073},
  the following theorem holds.

\vspace{3mm}
\bt[= Theorem A] \label{z078}\leavevmode\\
 \vspace{-5mm}
 \begin{itemize}
\item[(i)]
 For each  positive integer $n$ with $1 < n <p-1$,
 there exists a canonical isomorphism
 \begin{equation}
 \Theta_{\mfs \mfl_n, \mfX_{/S}}^\bigstar : \mfO \mfp^\ZZZ_{\mfs \mfl_n, \mfX_{/S}} \isom \mfO \mfp^\ZZZ_{\mfs \mfl_{(p-n)}, \mfX_{/S}} 
 \end{equation}
(i.e., the first isomorphism in (\ref{w4662}))
over $S$ satisfying that $\Theta_{\mfs \mfl_{(p-n)}, \mfX_{/S}}^\bigstar \circ   \Theta_{\mfs \mfl_n, \mfX_{/S}}^\bigstar = \mr{id}$.
\item[(ii)]
If, moreover, $r >0$ and we are given an element $\vec{\rho} \in (\mfc_n(\mbF_p)^\circledast)^{\times r}$, then we obtain, by restricting  $\Theta_{\mfs \mfl_n, \mfX_{/S}}^\bigstar$,  a canonical  isomorphism 
\begin{equation}
 \Theta_{\mfs \mfl_n, \vec{\rho}, \mfX_{/S}}^\bigstar : \mfO \mfp^\ZZZ_{\mfs \mfl_n, \vec{\rho}, \mfX_{/S}} \isom \mfO \mfp^\ZZZ_{\mfs \mfl_{(p-n)}, \vec{\rho}^{\, \bigstar}, \mfX_{/S}} 
\end{equation}
(i.e., the second  isomorphism in (\ref{w4662}))
over $S$ satisfying that $\Theta_{\mfs \mfl_{(p-n)}, \vec{\rho},\mfX_{/S}}^\bigstar \circ   \Theta_{\mfs \mfl_n, \vec{\rho}^{\, \bigstar},\mfX_{/S}}^\bigstar = \mr{id}$.
\end{itemize}
\et
\vspace{3mm}

Also, by Corollary \ref{z073ffff},  the following  assertion holds.

\vspace{3mm}
\bt[= Theorem B] \label{z07jjj8}\leavevmode\\
 \ \ \ The structure morphism $\mfO \mfp^\ZZZ_{\mfs \mfl_{(p-1)}, \mfX_{/S}} \migi S$ of $\mfO \mfp^\ZZZ_{\mfs \mfl_{(p-1)}, \mfX_{/S}}$ is an isomorphism.
That is to say, there exists a unique (up to isomorphism) dormant $\mfs \mfl_{(p-1)}$-oper on $\mfX_{/S}$.
 \et

\vspace{5mm}
\subsection{} \label{z0762}

In order to achieve
 a detailed understanding of the moduli stack of $\mfs \mfl_{(p-n)}$-opers (of the case where $n$ is, in a certain sense,  sufficiently  small relative to $p$),
 we will use, in this subsection, 
  Theorem \ref{z078}  and some  results concerning the moduli stack of dormant  $\mfs \mfl_{n}$-opers   obtained in ~\cite{Wak5}.

First, we shall  recall  the theory of dormant operatic fusion rings $\mfF^\ZZZ_{p, \mfs \mfl_n}$ discussed in ~\cite{Wak5}, \S\,7.
For each integer $n$ with $1 < n < p-1$ and $\vec{\rho} := (\rho_i)_{i=1}^r\in \mfc_n (\mbF_p)^{\times r}$ (where $\vec{\rho}:= \emptyset$ if $r =0$), write 
\begin{equation}
\mfO \mfp^\ZZZ_{\mfs \mfl_n, \vec{\rho}, g,r}
\end{equation}
for $\mfO \mfp^\ZZZ_{\mfs \mfl_n, \vec{\rho}, \mfX_{/S}}$ of the case where the pointed stable curve $\mfX_{/S}$ is taken to be the tautological pointed stable curve 
\begin{equation}
(f_{\mft \mfa \mfu} : \mfC_{g,r} \migi \overline{\mfM}_{g,r}, \{ \mfs_i,  :   \overline{\mfM}_{g,r} \migi \mfC_{g,r}\}_{i=1}^r)
\end{equation}
 (cf. \S\,\ref{z04}) over $\overline{\mfM}_{g,r}$.
That is,  $\mfO \mfp^\ZZZ_{\mfs \mfl_n, \vec{\rho}, g,r}$ is defined to be the stack in groupoids over $\mfS \mfc \mfh_{/\mr{Spec}(k)}$ whose category of sections over a $k$-scheme $S$ is the groupoid of the pairs $(\mfX_{/S}, \mcE^\spadesuit)$ consisting of a pointed stable curve $\mfX_{/S}$ over $S$ of type $(g,r)$ and a dormant  $\mfs \mfl_n$-opers on $\mfX_{/S}$ of radii $\vec{\rho}$.
(Indeed, it follows from ~\cite{Wak5}, Proposition 2.2.5, that any (dormant) $\mfs \mfl_n$-oper does not have  nontrivial automorphisms.)

According to  ~\cite{Wak5}, Theorem C and Theorem G,  $\mfO \mfp^\ZZZ_{\mfs \mfl_n, \vec{\rho}, g,r}$ is finite and generically \'{e}tale over  $\overline{\mfM}_{g,r}$.
Let   
\begin{equation}
N^\ZZZ_{p, \mfs \mfl_n, \vec{\rho}, g,r}
\end{equation}
 denote the generic degree of $\mfO \mfp^\ZZZ_{\mfs \mfl_n, \vec{\rho}, g,r}$ over $\overline{\mfM}_{g,r}$.
For each finite set $I$, $\mbN^I$
 denotes the free commutative monoid generated by $I$,
 and moreover, for each integer $l$,
 $\mbN^I_{\geq l}$
  denotes the submonoid of  $\mbN^I$ consisting of elements
$x = \sum_{i =1}^m a_i \lambda_i $ (where $\lambda_i \in I$ and $a_i \in \mbN$ for each $i = 1, \cdots, m$) with  $\sum_{i=1}^m  a_i \geq l$.
Here, 
recall from ~\cite{Wak5}, \S\,5.8, that
there exists an involution $\lambda \mapsto \lambda^\veebar$ on $\mfc_n (\mbF_p)$ that comes, via $\mfs \mfl_n \isom \mfp \mfg \mfl_n \stackrel{\chi_n}{\migi} \mfc_n$, from the involution on $\mfs \mfl_n$ given by assigning $(a_{i, j})_{i, j} \mapsto (- a_{n+1 -j, n +1 -i})_{i, j}$.
This involution extends by linearity to an involution $x \mapsto x^\veebar$ on $\mbN^I_{\geq l}$ (for each $l \in \mbZ$).
Then, the function
\begin{align} \label{e0e0e}
N^\ZZZ_{p, \mfs \mfl_n, g} : \ \  \mbN_{\geq 3-2g}^{\mfc_{n} (\mbF_p)} &\migi \ \  \mbZ \\
\sum_{i=1}^r\rho_i  \ \ & \mapsto N^\ZZZ_{p, \mfs \mfl_n, \vec{\rho}, g,r} \notag
\end{align}
is verified to be well-defined, i.e., the value $N^\ZZZ_{p, \mfs \mfl_n, \vec{\rho}, g,r}$ does not depend on the ordering of $\rho_1, \cdots, \rho_r$ (cf. ~\cite{Wak5}, the discussion following Proposition 7.5.2).
It follows from ~\cite{Wak5},  Theorem 7.10.4 (i), that 
if $g_1$, $g_2$ are nonnegative integers and $x \in \mbN_{\geq 3-2g_1}^{\mfc_{n} (\mbF_p)}$, $y \in \mbN_{\geq 3-2g_2}^{\mfc_{n} (\mbF_p)}$, then the collection of functions
$\{ N^\ZZZ_{p, \mfs \mfl_n, g}\}_{g \geq 0}$ satisfies the following rule:
\begin{equation}
N^\ZZZ_{p, \mfs \mfl_n, g_1 +g_2} (x +y) = \sum_{\lambda \in \mfc_n (\mbF_p)} N^\ZZZ_{p, \mfs \mfl_n, g_1} (x + \lambda) \cdot N^\ZZZ_{p, \mfs \mfl_n, g_2} (y + \lambda^\veebar).
\end{equation}
Also,
 it follows from ~\cite{Wak5}, Theorem E, that
for any nonnegative integer $g$ and any $x \in \mbN_{\geq 3-2g}^{\mfc_{n} (\mbF_p)}$,
the following equality holds:
\begin{align}
N^\ZZZ_{p, \mfs \mfl_n, g} (x^\veebar) =  N^\ZZZ_{p, \mfs \mfl_n, g} (x).
\end{align}
In particular, 
the functions $N^\ZZZ_{p, \mfs \mfl_n, 0}$ 
forms a pseudo-fusion rule (cf. ~\cite{Wak5}, Definition 7.6.1) on the finite set $\mfc_{n} (\mbF_p)$ (with  involution $x \mapsto x^\veebar$), and hence, 
one obtains the fusion ring 
\begin{equation} \label{d0d0d}
\mfF^\ZZZ_{p, \mfs \mfl_n}
\end{equation}
associated with $N^\ZZZ_{p, \mfs \mfl_n, 0}$ (cf. ~\cite{Wak}, (834)).

The decomposition (\ref{w4660}) implies that 
$\mfc (\mbF_p) \setminus \mfc_n (\mbF_p)^\circledast$ is contained in the kernel (cf. ~\cite{Wak5}, Remark 7.7.1) of $N^\ZZZ_{p, \mfs \mfl_n, 0}$, i.e., $N^\ZZZ_{p, \mfs \mfl_n, \vec{\rho}, 0,r} \neq 0$ only if $\vec{\rho} \in (\mfc_n(\mbF_p)^\circledast)^{\times r}$.
Hence, the restriction 
\begin{equation} \label{b00i9}
N'_n := N^\ZZZ_{p, \mfs \mfl_n, 0} |_{\mbN_{\geq 3}^{\mfc_n (\mbF_p)^\circledast}}
\end{equation}
of $N^\ZZZ_{p, \mfs \mfl_n, 0}$ to $\mbN_{\geq 3}^{\mfc_n (\mbF_p)^\circledast} \subseteq \mbN_{\geq 3}^{\mfc_n (\mbF_p)}$ forms a pseudo-fusion rule on $\mfc_n (\mbF_p)^\circledast$, and the natural inclusion $\mfF_{N'_n} \migiincl \mfF^\ZZZ_{p, \mfs \mfl_n}$ gives rise to an isomorphism
\begin{equation} \label{w4881}
(\mfF_{N'_n })_{\mr{red}} \isom (\mfF^\ZZZ_{p, \mfs \mfl_n})_{\mr{red}}
\end{equation}
(cf. ~\cite{Wak5}, Remark 7.7.1) between the reduced  rings associated with $\mfF_{N'_n}$ and $\mfF^\ZZZ_{p, \mfs \mfl_n}$ respectively.
Let
\begin{equation}
\overline{\mcC as}_{p, \mfs \mfl_n}
\end{equation}
be  the  element of $ (\mfF^\ZZZ_{p, \mfs \mfl_n})_{\mr{red}}$ defined to be the image of $\sum_{\lambda \in \mfc_n (\mbF_p)} \lambda^\veebar \cdot \lambda \in \mfF^\ZZZ_{p, \mfs \mfl_n}$ via the quotient $ \mfF^\ZZZ_{p, \mfs \mfl_n} \migisurj (\mfF^\ZZZ_{p, \mfs \mfl_n})_{\mr{red}}$.
Then, by means of  Theorem  \ref{z078}, 
we have the following  corollary.
In particular, one may extend (cf. Corollary \ref{z079} (ii)) a result in the paper ~\cite{Wak5}  (cf.  ~\cite{Wak5}, Theorem H)
concerning an explicit computation of  the value  $N^\ZZZ_{p, \mfs \mfl_n, \emptyset, g,0}$.

\vspace{3mm}
\bco \label{z079}\leavevmode\\
\vspace{-5mm}
\begin{itemize}
\item[(i)]
There exists a canonical isomorphism
\begin{equation}
(\mfF^\ZZZ_{p, \mfs \mfl_n})_{\mr{red}} \isom (\mfF^\ZZZ_{p, \mfs \mfl_{(p-n)}})_{\mr{red}} \ \left(=:\mfF\right)
\end{equation}
of rings that sends $\overline{\mcC as}_{p, \mfs \mfl_n}$ to $\overline{\mcC as}_{p, \mfs \mfl_{(p-n)}}$ ($=: \overline{\mcC as}$).
In particular, for each $\vec{\rho} = (\rho_i)_{i=1}^r \in (\mfc_n (\mbF_p)^\circledast)^{\times r}$ (where $\vec{\rho}:=\emptyset$ if $r =0$),  the following equalities hold:
\begin{equation}
N_{p, \mfs \mfl_n, \vec{\rho}, g, r}^\ZZZ  = N_{p, \mfs \mfl_{(p-n)}, \vec{\rho}^{\, \bigstar}, g, r}^\ZZZ   = \sum_{\chi \in \mr{Hom}(\mfF, \mbC)} \chi (\overline{\mcC as})^{g-1}\cdot \prod_{i=1}^r \chi(\rho_i),
\end{equation}
where  $\mr{Hom} (\mfF, \mbC)$ denotes the set of ring homomorphisms $\mfF \migi \mbC$.

\item[(ii)]
Suppose that $p>n \cdot \mr{max} \{g-1, 2 \}$.
  Then, the generic degrees
  $N^\ZZZ_{p, \mfs \mfl_n, \emptyset, g, 0}$ and   $N^\ZZZ_{p, \mfs \mfl_{(p-n)}, \emptyset, g, 0}$
 are given by the following formula:
 \vspace{3mm}
 \begin{equation} 
N^\ZZZ_{p, \mfs \mfl_n, \emptyset, g, 0} =  N^\ZZZ_{p, \mfs \mfl_{(p-n)}, \emptyset, g, 0}   =    \frac{ p^{(n-1)(g-1)-1}}{n!} \cdot  
 \sum_{\genfrac{.}{.}{0pt}{}{(\zeta_1, \cdots, \zeta_n) \in \mbC^{\times n} }{ \zeta_i^p=1, \ \zeta_i \neq \zeta_j (i\neq j)}}
\frac{(\prod_{i=1}^n\zeta_i)^{(n-1)(g-1)}}{\prod_{i\neq j}(\zeta_i -\zeta_j)^{g-1}}.  
 \end{equation}

\end{itemize}
\eco
\begin{proof}
Assertion (i) follows from  ~\cite{Wak5} Theorem F,  Theorem \ref{z078} (ii), the definitions of $N^\ZZZ_{p, \mfs \mfl_n, 0}$ and $\mfF^\ZZZ_{p, \mfs \mfl_n}$, and the isomorphism (\ref{w4881}) (as well as the isomorphism (\ref{w4881}) of the case where the integer ``$n$" is taken to be ``$p-n$").
Assertion (ii) follows from ~\cite{Wak5}, Theorem H, and    Theorem \ref{z078} (i).
\end{proof}

\vspace{5mm}
\subsection{} \label{z0762}

In the following, we focus on the case where $n = p-2$.
By means of  Theorem \ref{z078} and results in the $p$-adic Teichm\"{u}ller developed by S. Mochizuki (cf. ~\cite{Mzk2}),  
one may  prove Corollaries  \ref{z080} and  \ref{z08022} described below.

\vspace{3mm}
\bco \label{z080}\leavevmode\\
 \ \ \ 
Let $\vec{\rho} \in \mfc_{(p-2)} (\mbF_p)^{\times r}$.
Then, the stack $\mfO \mfp^\ZZZ_{\mfs \mfl_{(p-2)}, \vec{\rho}, g,r}$ is a (possibly empty) geometrically connected, proper, and smooth Deligne-Mumford stack over $k$ of dimension $3g-3+r$.
Moreover, the natural morphism $\mfO \mfp^\ZZZ_{\mfs \mfl_{(p-2)}, \vec{\rho}, g,r} \migi \overline{\mfM}_{g,r}$ is finite, faithfully flat, and generically \'{e}tale. 
\eco
\begin{proof}
The assertion follows from Theorem \ref{z078} (ii) of the case $n =2$ and  ~\cite{Mzk2}, Chap.\,II, \S\,2.8, Theorem 2.8, which asserts that 
 $\mfO \mfp^\ZZZ_{\mfs \mfl_2, \vec{\rho}^{\, \bigstar}, g,r}$ 
 satisfies the same properties as the properties  desired  for $\mfO \mfp^\ZZZ_{\mfs \mfl_{(p-2)}, \vec{\rho}, g,r}$ in the statement.
\end{proof}

\vspace{3mm}

Next, we consider the structure of the reduced ring $(\mfF^\ZZZ_{p, \mfs \mfl_{(p-2)}})_{\mr{red}}$
 associated with $\mfF^\ZZZ_{p, \mfs \mfl_{(p-2)}}$.
(In order to perform any computation that we will need in the ring $\mfF^\ZZZ_{p, \mfs \mfl_{(p-2)}}$, it suffices to understand the structure of $(\mfF^\ZZZ_{p, \mfs \mfl_{(p-2)}})_{\mr{red}}$ (cf. Corollary  \ref{z079} (i))).

Let us write
\begin{equation}
\mbF := \{ a \in \mbZ \ | \ 0 \leq a \leq \frac{p-3}{2} \}.
\end{equation}
The composite
\begin{equation} \label{b00t67}
\mbF \stackrel{}{\migi} 2^{\mbF_p}_{\sharp (-) = (p-2)} \migi \mfc_{(p-2)} (\mbF_p)
\end{equation}
is injective, where the first arrow denotes the map 
\begin{align}
\mbF &\migi 2^{\mbF_p}_{\sharp (-) = (p-2)}  \\
a &\mapsto  \mbF_p \setminus \{ 0, \overline{-2a -1} \}. \notag
\end{align}
We shall regard $\mbF$ as a subset of $\mfc_{(p-2)} (\mbF_p)$ via this composite injection.
\vspace{3mm}
\bco \label{z08022}\leavevmode\\
 \ \ \ 
The complement of $\mbF$ in $\mfc_{(p-2)} (\mbF_p)$ is contained in the kernel (cf. ~\cite{Wak5}, Remark 7.7.1) of $N^\ZZZ_{p, \mfs \mfl_{(p-2)}, 0}$.
Moreover, the structure constants $N^\ZZZ_{p, \mfs \mfl_{(p-2)}, 0} (\alpha + \beta + \gamma)$ (where $\alpha$, $\beta$, $\gamma \in \mfc_{(p-2)} (\mbF_p)$) of the fusion ring $\mfF^\ZZZ_{p, \mfs \mfl_{(p-2)}}$ are given as follows:
 \begin{equation}
N^\ZZZ_{p, \mfs \mfl_{(p-2)}, 0} (\alpha + \beta + \gamma) = \begin{cases} 1 & \text{if $(\alpha, \beta,  \gamma) \in W$,} \\
0 & \text{if $(\alpha, \beta,  \gamma) \notin W$,}
\end{cases}
\end{equation}
where $W$ denotes a subset  of $\mbF^{\times 3}$ ($\subseteq \mfc_{(p-2)} (\mbF_p)^{\times 3}$) defined  to be
 \begin{equation}
 W := \Bigg\{ (s, t, u) \in \mbF^{\times 3} \   \Biggm| 
  \text{  $s+ t+u \leq p-2$,  and } \   \begin{matrix}
 0 \leq s \leq t + u,  \\
 0 \leq t \leq u+s,  \\ 
 0 \leq u\leq s +t.
 \end{matrix}
   \Biggm\}.
 \end{equation}
\eco
\begin{proof}

Consider the composite
\begin{equation}
\mbF \migi  \mfs \mfl_2 (\mbF_p) \stackrel{\chi_2 (\mbF_p)}{\migi} \mfc_2 (\mbF_p),
\end{equation}
where the first  arrow denotes the map given by assigning $a \mapsto \begin{pmatrix} \overline{a} & 0 \\ 0 & -\overline{a}\end{pmatrix}$.
This composite factors through the inclusion $\mfc_2 (\mbF_p)^\circledast \migi \mfc_2 (\mbF_p)$.
By passing to the resulting injection $w : \mbF \migiincl \mfc_2 (\mbF_p)^\circledast$,  we regard $\mbF$ as a subset of $\mfc_2 (\mbF_p)^\circledast$.
According to the discussion in  ~\cite{Wak5}, \S\,7.11 (or, ~\cite{Mzk2}, Introduction, \S\,1.2, Theorem 1.3),  the structure constants $N^\ZZZ_{p, \mfs \mfl_{2}, 0} (\alpha + \beta + \gamma) $  (where $\alpha$, $\beta$, $\gamma \in \mfc_2  (\mbF_p)^\circledast$) of $\mfF_{N'_2}$ (cf. (\ref{b00i9}))  are given as follows:
 \begin{equation}
N^\ZZZ_{p, \mfs \mfl_{2}, 0} (\alpha + \beta + \gamma)  = \begin{cases} 1 & \text{if $(\alpha, \beta,  \gamma)  \in W$,} \\
0 & \text{if $(\alpha, \beta,  \gamma)  \notin W$.}
\end{cases}
\end{equation}
On the other hand, the composite 
\begin{equation}
\mbF \stackrel{w}{\migiincl} \mfc_2(\mbF_p)^\circledast  \stackrel{\theta_2^\bigstar}{\migi} \mfc_{(p-2)} (\mbF_p)^\circledast  \migiincl \mfc_{(p-2)} (\mbF_p)
\end{equation}
coincides with the composite (\ref{b00t67}).
Thus, the assertion follows from Corollary  \ref{z079} (i) of the case $n =2$.
\end{proof}

\end{document}